\documentclass[11pt]{article}
\usepackage{epsfig,amsmath,latexsym}
\usepackage{amsfonts}

\newtheorem{theorem}{Theorem}[section]

\newtheorem{defi}{Definition}[section]
\newtheorem{lemma}[theorem]{Lemma}

\catcode`\@=11
\renewcommand{\section}{
        \setcounter{equation}{0}
        \@startsection {section}{1}{\z@}{-3.5ex plus -1ex minus
        -.2ex}{2.3ex plus .2ex}{\large\bf}
        }
\catcode`@=12

\def\slfrac#1#2{\hbox{\kern.1em %
 \raise.5ex\hbox{\the\scriptfont0 #1}\kern-.11em %
 /\kern-.15em\lower.25ex\hbox{\the\scriptfont0 #2}}}

\newcommand{\pf}{\noindent{\bf Proof.~}}
\newcommand{\beq}{\begin{eqnarray}}
\newcommand{\eeq}{\end{eqnarray}}
\newcommand{\beql}[1]{\begin{eqnarray}\label{#1}}
\newcommand{\beqs}{\begin{eqnarray*}}
\newcommand{\eeqs}{\end{eqnarray*}}
\newcommand{\eqn}[1]{(\ref{#1})}

\newcommand{\lf}{\lfloor}
\newcommand{\rf}{\rfloor}

\newcommand{\cc}{{\mathbb C}}

\newcommand{\rr}{{\mathbb R}}
\newcommand{\zz}{{\mathbb Z}}

\newcommand{\nn}{{\mathbb N}}

\newcommand{\bv}{{\mathbf v}}
\newcommand{\bw}{{\mathbf w}}

\newcommand{\sA}{{\mathcal A}}

\newcommand{\sD}{{\mathcal D}}

\newcommand{\sL}{{\mathcal L}}
\newcommand{\sP}{{\mathcal P}}
\makeatletter
\def\section{\@startsection {section}{1}{\z@}{-3.5ex plus -1ex minus 
 -.2ex}{2.3ex plus .2ex}{\normalsize\bf}}
\def\subsection{\@startsection {subsection}{1}{\z@}{-3.5ex plus -1ex minus
 -.2ex}{2.3ex plus .2ex}{\normalsize\bf}}
\def\@sect#1#2#3#4#5#6[#7]#8{\ifnum #2>\c@secnumdepth
     \def\@svsec{}\else
     \refstepcounter{#1}\edef\@svsec{\csname the#1\endcsname.\hskip .75em }\fi
     \@tempskipa #5\relax
      \ifdim \@tempskipa>\z@
        \begingroup #6\relax
          \@hangfrom{\hskip #3\relax\@svsec}{\interlinepenalty \@M #8\par}%
        \endgroup
       \csname #1mark\endcsname{#7}\addcontentsline
         {toc}{#1}{\ifnum #2>\c@secnumdepth \else
                      \protect\numberline{\csname the#1\endcsname}\fi
                    #7}\else
        \def\@svsechd{#6\hskip #3\@svsec #8\csname #1mark\endcsname
                      {#7}\addcontentsline
                           {toc}{#1}{\ifnum #2>\c@secnumdepth \else
                             \protect\numberline{\csname the#1\endcsname}\fi
                       #7}}\fi
     \@xsect{#5}}
\def\@begintheorem#1#2{\it \trivlist \item[\hskip \labelsep{\bf #1\ #2.}]}
\makeatother

\begin{document}

\begin{center}
{\Large {\bf Apollonian Circle Packings: Number Theory}} \\

\vspace{1.5\baselineskip}
{\em Ronald L. Graham}~\footnote{Current address: Dept. of Computer Science,
Univ. of Calif. at San Diego, La Jolla, CA 92093} \\
\vspace*{.2\baselineskip}
{\em Jeffrey C. Lagarias}~\footnote{Work partly done during a visit to
the Institute for Advanced Study.} \\
\vspace*{.2\baselineskip}
{\em Colin L. Mallows} \\
\vspace*{.2\baselineskip}
{\em Allan R. Wilks} \\
\vspace*{.2\baselineskip}
AT\&T Labs, 
Florham Park, NJ 07932-0971 \\
\vspace*{1.5\baselineskip}

{\em Catherine H. Yan} \\
\vspace*{.2\baselineskip}
Texas A\&M University, 
College Station, TX 77843\\
\vspace*{1.5\baselineskip}
%\today
(November 1, 2002 version)

\vspace{1.5\baselineskip}
{\bf ABSTRACT}
\end{center}

Apollonian circle packings arise by repeatedly filling the interstices
between mutually tangent circles with further tangent circles.  It is
possible for every circle in such a packing to have integer radius of
curvature, and we call such a packing an {\em integral Apollonian 
circle packing.} This paper studies number-theoretic properties of
the set of integer curvatures appearing in such packings. 
Each Descartes quadruple of four tangent circles in the
packing gives an integer
solution to the Descartes equation, which relates the 
radii of curvature of four
mutually tangent circles:  
$ x^2 + y^2 + z^2 + w^2 = \frac{1}{2}(x + y + z +w)^2.$ 
Each integral Apollonian circle packing is classified by a
certain {\em root quadruple} of integers that satisfies the Descartes
equation, and that corresponds to a particular quadruple of circles
appearing in the packing. 
We express the number of root quadruples with fixed minimal
element $-n$ as a class number, and give an exact formula for it.
 We study which
integers occur in a given integer packing, and determine congruence
restrictions which sometimes apply. 
We present evidence suggesting that
the set of integer radii of curvatures that
appear in an integral Apollonian circle packing has
positive density, and in fact represents all sufficiently large integers
not excluded by congruence conditions. Finally, we discuss
asymptotic properties of the set of curvatures obtained
as the packing is recursively constructed from a root quadruple.
%In related work,  ``Apollonian Packings: Geometry
%and Group Theory,'' we investigate a variety of group-theoretic
%properties of these configurations, as well as various extensions
%to higher dimensions and other spaces, such as hyperbolic space.

\vspace*{1.5\baselineskip}
\noindent
Keywords: Circle packings, Apollonian circles, Diophantine equations 
\newpage
\setcounter{page}{1}
\begin{center}
{\Large {\bf  Apollonian Circle Packings: Number Theory }} \\

\vspace{3\baselineskip}

\end{center}
\setlength{\baselineskip}{1.5\baselineskip}

%--------------------------------------------------------------------
%
% Section 1
%
%--------------------------------------------------------------------
%
\section{Introduction}
\begin{figure}[htbp]\label{std}
\centerline{\epsfxsize=4.0in \epsfbox{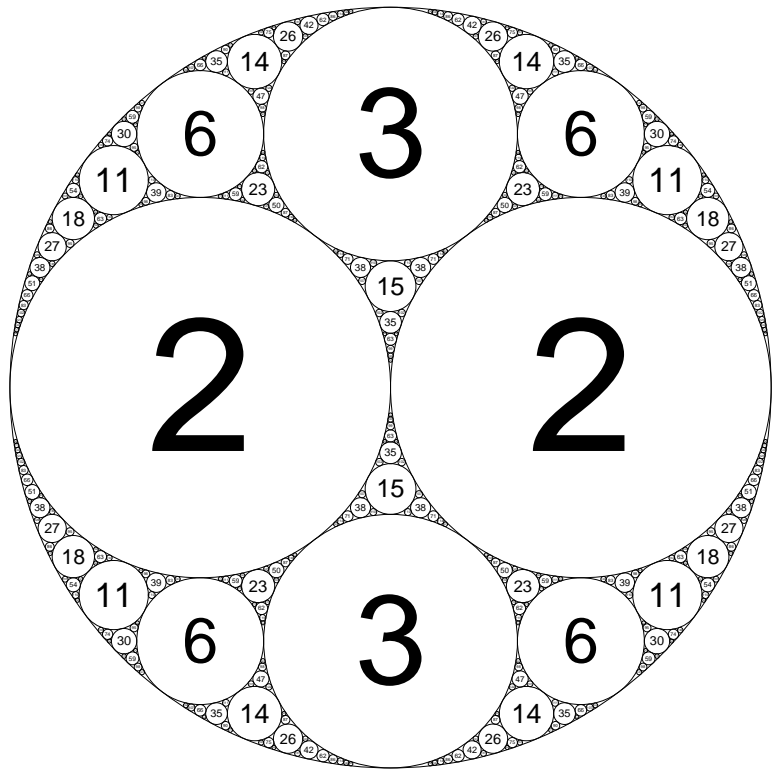}}
\caption{The integral Apollonian circle packing (-1,~2,~2,~3)}
\end{figure}

Place two tangent circles of radius 1/2 inside and tangent to a circle
of radius 1.  In the two resulting curvilinear triangles fit tangent
circles as large as possible.  Repeat this process for the six new
curvilinear triangles, and so on.  The result is pictured in
Figure \ref{std}, where each circle has been labeled with 
its curvature---the
reciprocal of its radius.

Remarkably, every circle in Figure \ref{std} has integer curvature.
Even more remarkable is that if the picture is centered at the origin
of the Euclidean plane with the centers of the ``2'' circles on the
$x$-axis, then each circle in the picture has the property that the
coordinates of its center, multiplied by its curvature, are also
integers. In this paper we are concerned with circle packings having
the first of these properties; the  latter property 
 is addressed in a companion paper \cite[Section 3]{GLMWY22}.

An {\em Apollonian circle packing} is any packing of circles constructed
recursively from an initial configuration of four mutually
tangent circles by the procedure above. More precisely, one
starts from a {\em Descartes configuration}, which is a set
of four mutually tangent circles with disjoint interiors, suitably
defined. In the example above, the enclosing circle has ``interior''
equal to its exterior, and its curvature is given a negative sign.
Recall that in a quadruple
of mutually touching circles the curvatures $(a,~b,~c,~d)$ satisfy the
{\em Descartes equation}
\beq\label{descartes}
a^2+b^2+c^2+d^2 = \frac{1}{2}(a+b+c+d)^2, 
\eeq
as observed by Descartes in 1643 (in an equivalent form).
Any quadruple $(a,~b,~c,~d)$ satisfying this equation is called
a {\em  Descartes quadruple.}
An {\em integral Apollonian circle packing} is an Apollonian circle
packing in which
every circle has an integer curvature. The starting point
of this paper is the observation that 
if an initial Descartes
configuration has all integral curvatures, then the whole
packing is integral, and conversely.
This integrality property of packings has been discovered
repeatedly; perhaps the first observation of it is in the 
1937 note of F. Soddy~\cite{Sod37}
``The bowl of integers and the Hexlet''. It is discussed in
some detail in Aharonov and Stephenson~\cite{AS97}.

In this paper we study
integral Apollonian circle packings viewed as equivalent under 
Euclidean motions,
an operation which preserves the curvatures of all circles.
Such packings are classified by their
root quadruple, a notion defined in \S3. 
This is the ``smallest'' quadruple in the packing as measured in
terms of curvatures of the circles. In the packing
above the root quadruple is $(-1,~2,~2, ~3),$ where $-1$ represents
the (negative) curvature of the bounding circle.
We study the set of integers (curvatures) represented by a packing 
using the {\em Apollonian group} $\sA$, which
is a subgroup of $GL(4, \zz)$ which
acts on integer Descartes
quadruples.  This action permits one to  ``walk around'' on a fixed
Apollonian packing, moving from one Descartes quadruple to any other
quadruple in the same packing. 
% as shown in \cite[Theorem 3.6]{GLMWY21}.
The Apollonian group  was introduced by
Hirst~\cite{Hirst67} in 1967, who used it bounding the Hausdorff dimension
of the residual set of an Apollonian packing; it was also used
in  S\"{o}derberg,\cite{So92} and Aharonov and Stephenson~\cite{AS97}.
Descartes quadruples associated to different root quadruples cannot
be reached by the action of $\sA$, and the action of the
Apollonian group partitions the set of integer Descartes quadruples
into infinitely many equivalence classes (according to which
integral Apollonian packing they belong.) 
 By scaling an integer
Apollonian packing by an appropriate homothety, one
may obtain a  {\em primitive
integral Apollonian packing,} which is one whose Descartes quadruples
have integer curvatures with greatest common divisor $1$. Thus 
the study of integral Apollonian packings essentially reduces to
the study of primitive packings.
 
The simplest integral Apollonian circle packing is the one with
root quadruple $(0,~0,~1,~1)$, which is pictured in Figure 2.
This packing is special in several ways.
It is degenerate in that it has two circles with
``center at infinity'', whose boundaries are straight lines,
and it is the only primitive integral Apollonian circle packing that
is unbounded. It is also the only primitive
integral Apollonian circle packing
that contains infinitely many copies of the root quadruple.
This particular packing has already played a role in number theory.
That part of the packing  
in an interval of length two between the tangencies of two
adjacent circles of radius one, consisting  of
the (infinite) set of circles tangent to one of the straight lines,
forms a set of  ``Ford circles'',
after shrinking all circles by a factor of two.  These
circles, introduced by Ford (see \cite{Fo38}),
 can be labelled by the Farey fractions on the interval
$(0,1)$ and used to prove basic results in one-dimensional
Diophantine approximation connected with the Markoff spectrum,
see Rademacher~\cite[pp. 41--46]{Ra64} and Nicholls~\cite{Ni78}.

\begin{figure}[htbp]\label{vertical}
\centerline{\epsfxsize=3.0in \epsfbox{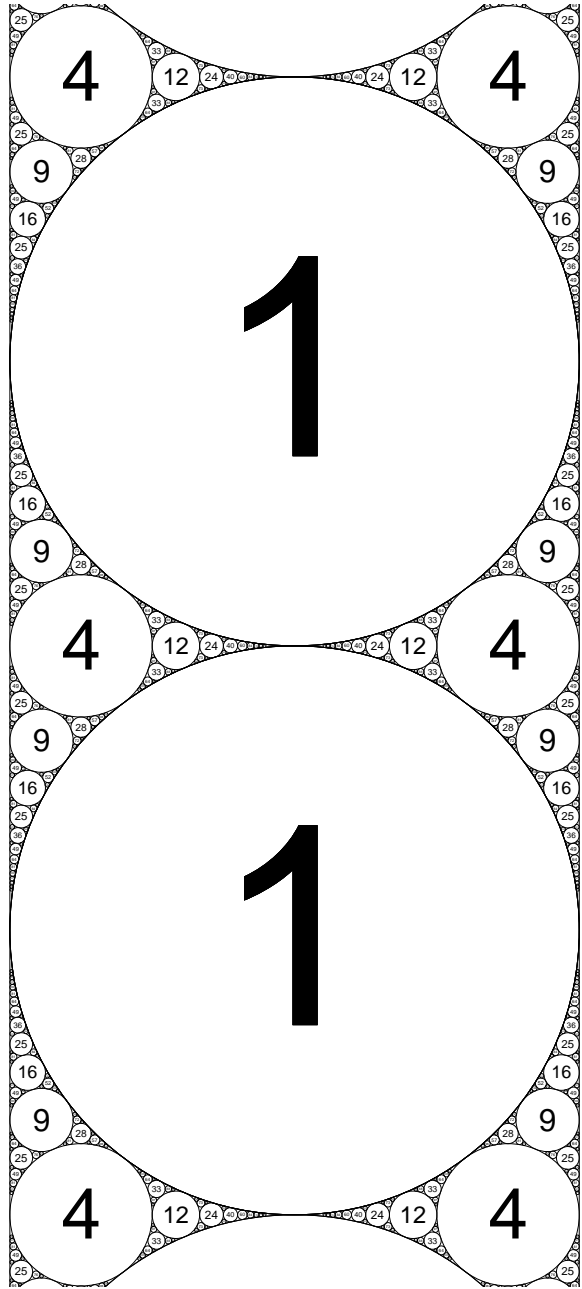}}
\caption{The Apollonian circle packing $(0,~0,~1,~1).$}
\end{figure}

In this paper our interest is in those
properties of the set of integral Apollonian circle
packings that are of a Diophantine nature. These include the 
distribution of integer Descartes quadruples, of integer
root quadruples, and
the representation and the distribution of the integers (curvatures) 
occurring in a fixed integral Apollonian circle packing.
Finally we consider the size
distribution of elements in the Apollonian group, a group of integer matrices
associated to such packings.

To begin with, the full set of all integer Descartes quadruples
(taken over all integral Apollonian packings)
is enumerated by the integer solutions to the Descartes equation
(up to a sign.) 
In \S2 we determine  asymptotics for the total number of integer solutions
to the Descartes equation of Euclidean norm below a given bound. 

In \S3 we define the Apollonian group.
We describe a reduction
theory which multiplies Descartes quadruples by elements
of this group and uses it to find a quadruple of
smallest size in a given packing, called a {\em root quadruple.}
We prove the existence and uniqueness
of a root quadruple associated to each
integral Apollonian packing.

In \S4 we study the root quadruples of primitive integer packings.
We show that the number $N_{root}^*(-n)$  of primitive root quadruples 
containing a given negative integer $-n$ 
is equal to the  class number of primitive integral binary quadratic forms
of discriminant $-4n^2$, under $GL(2, \zz)$-equivalence. Using
this fact we obtain in Theorem~\ref{Nth42}
an exact formula for the number of such quadruples,
which is
\beql{105}
 N_{root}^*(-n)= \frac{n}{4} \prod_{p|n}(1-\frac{\chi_{-4}(p)}{p})+
2^{\omega(n)-\delta_n-1},
\eeq
where $\chi_{-4}(n) = (-1)^{(n-1)/2}$ for $n$ and $0$ for even $n$,
$\omega(n)$ denotes  the number of distinct primes dividing $n$, and 
$\delta_n =1$ if $n \equiv 2(\bmod~4)$ and $\delta_n=0$
otherwise.
This formula was independently discovered and
proved by S. Northshield~\cite{No02}. 
As a consequence we show that   $N_{root}^*(-n)$  has an  upper
bound $O( n \log \log n)$ and a lower bound  
$\Omega (\frac {n}{(\log \log n)}),$ respectively. These bounds
are sharp up to multiplicative constants.

In \S5 we study the integer curvatures appearing in
a single integral Apollonian packing, counting integers with
multiplicity.
D. Boyd~\cite{Bo82} showed that the number of circles occurring
in a bounded Apollonian packing having curvature less than
a bound $T$ grows like $T^{\alpha + o(1)}$, where
$\alpha \approx 1.30...$ is the Hausdorff dimension of the residual
set of any  Apollonian circle packing. This result applies
to integral Apollonian packings. We observe that these
 integers can 
be put in one-to-one correspondence with elements of the Apollonian
group, using the root quadruple. 
Using this result
we show that the number of elements of the Apollonian group which
have norm less than $T$ is of order $T^{\alpha + o(1)}$,
as $T \to \infty.$ 
 
In \S6 we study the integer curvatures appearing in a packing,
counted without multiplicity. 
We show that there are always
nontrivial congruence restrictions $(mod~24)$ on the integers
that occur.
We give some evidence suggesting
that such congruence restrictions can only
involve powers of the primes $2$ and $3$.
We conjecture that in any integral Apollonian packing, all sufficiently large
integers occur, provided they are not excluded by a congruence condition. 
This may be a hard problem, however, since we show that it is analogous
to Zaremba's conjecture stating that there is a fixed integer $K$
such that for all denominators $n \geq 2$ there is a rational $\frac{a}{n}$
in lowest terms whose continued fraction expansion has all partial
quotients bounded by $K$.

In \S7 we study the set 
of integer curvatures at ``depth $n$'' in an
integral Apollonian packing, where $n$ measures the distance to the root
quadruple. There are exactly
$4 \times 3^{n - 1}$ such elements. We determine the maximal and minimal
curvature in this set, and also formulate
 a conjecture concerning the asymptotic
size of the median curvature as $ n \to \infty.$
These problems are  related to the joint spectral
radius of the matrix generators  
$\Sigma = \{ S_1, S_2, S_3, S_4\} $ of the Apollonian group, which
we determine.

In \S8 we conclude with some directions for
further work and  a list of open problems.

There has been extensive previous work on various aspects of Apollonian
circle packings, related to geometry, group theory and fractals.
The name ``Apollonian packing'' traces back
to Kasner and Supnick~\cite{KS43} in 1943.
It has been popularized  by Mandlebrot~\cite[p. 169ff]{M82},
who observed a connection with work of Apollonius of Perga,
around 200BC, see also Coxeter~\cite{Cox68}. 
 Further discussion and references can be found in 
 Aharonov and Stephenson~\cite{AS97} and Wilker~\cite{Wi81}.
See also the companion papers ~\cite{GLMWY21}, \cite{GLMWY22},\cite{GLMWY23}
and \cite{LMW00}, where we investigate a variety of group-theoretic
properties of these configurations, as well as various extensions
to higher dimensions and other spaces, such as hyperbolic space.

\paragraph{Notation.} We use several different measures
of the size of a vector $\bw= (w_1, w_2, w_3, w_4)$.
In particular we set $L(\bw):= w_1+w_2+w_3+w_4$, a
measure of size which is not a norm. We let 
$H(\bw) := (w_1^2 + w_2^2 + w_3^2 +w_4^2)^{1/2}$ or $|\bw|$ denote
the Euclidean norm, while
$||\bw||_\infty := \max_{1 \le i \le 4} |w_i|$ 
denotes the supremum norm of $\bw$.
%--------------------------------------------------------------------
%
% Section 2
%
%--------------------------------------------------------------------
%

\section{Integral Descartes Quadruples}
\setcounter{equation}{0}
An Apollonian circle  packing is {\em integral} if 
every circle of the packing has an integer curvature. {From} 
\eqref{descartes} it follows that if $a, b, c,$ are given, the
curvatures $d, d'$ of the two circles that are tangent to all three
satisfy
\beqs
d, \ d'=a+b+c\pm 2q_{abc}, 
\eeqs
where 
\beqs 
q_{abc}=\sqrt{ab+bc+ac}. 
\eeqs
Hence 
\beq\label{duplex}
d+d'=2(a+b+c).  
\eeq
In other words, given four mutually tangent circles
with curvatures $a, b, c,d$, the
curvature of the other circle that touches the first three is given by
\beq\label{stepa}
d'=2a+2b+2c-d.
\eeq
It follows that an Apollonian packing is integral if the starting Descartes 
quadruple consists entirely of integers.

The relation (\ref{duplex}) is the basis of the integrality property
of Apollonian packings. 
It generalizes to $n$ dimensions, where the curvatures $X_i$ of a
set of $n+1$ mutually tangent spheres in $\rr^n$ (having distinct
tangents) are related to the
curvatures $X_0$ and $X_{n+2}$ of the two spheres that are tangent
of all of these by
$$
X_0 + X_{n+2} = \frac{2}{n - 1}( X_1 + X_2 + \cdots + X_n).
$$
This relation gives integrality in dimensions $n = 2$ and $n = 3$;
the three dimensional case is studied in Boyd~\cite{Bo73}.
It even generalizes further to sets of equally inclined spheres with
inclination parameter $\gamma$, with the constant
$\frac{2}{n + \frac{1}{\gamma}}$; the case $\gamma = -1$ is the
mutually tangent case, cf. Mauldon~\cite{Mau62} and 
Weiss~\cite[Theorem 3]{We83}.

\begin{defi}
{\em 

(i) An {\em integer Descartes quadruple} 
${\bf a} = (a, b, c, d) \in \zz^4$ is any
integer representation of zero by the indefinite integral quaternary
quadratic form, 
$$ Q_{\sD}(w, x, y, z) := 2(w^2 + x^2 + y^2 + z^2 ) - (w  + x + y + z)^2, $$
which we call the {\em Descartes integral form}. That is, writing
$\bv = (w,~x,~y,~z)^T,$ we have $ Q_{\sD}(w, x, y, z)= \bv^T Q_{\sD} \bv,$
where
\beql{300a}
Q_{\sD} = \left[ \begin{array}{cccc}
           ~1 &  -1 &  -1 &  -1 \\
          -1 &   ~1 &  -1 &  -1 \\
          -1 &  -1 &   ~1 &  -1 \\
          -1 &  -1 &  -1 &   ~1\end{array} \right].
\eeq
This quadratic form
has determinant $det (Q_{\sD}) = - 16$ and, on identifying the
form $Q_{\sD}$ with its symmetric integral matrix, it satisfies 
$Q_{\sD}^2 = 4I.$

(ii) An integer Descartes quadruple is {\em primitive} if
$gcd(a, b, c, d) = 1.$
}
\end{defi}

\noindent In studying the geometry of 
Apollonian packings 
(\cite{GLMWY21}-\cite{GLMWY23}, \cite{LMW00}) we use instead
 a scaled version of the Descartes integral form, namely
the {\em Descartes quadratic form} $Q_2 := \frac{1}{2} Q_{\sD}.$

\begin{defi}
{\em 
The size of any real quadruple $(a, b, c, d) \in \rr^4$
is measured by the
 {\em Euclidean  height} $H({\bf a})$, which is:
\beql{301}
H({\bf a}) := (a^2 + b^2 + c^2 + d^2)^{1/2}.
\eeq
}
\end{defi}

Now let $N_{\sD}(T)$ count the number of integer Descartes quadruples 
with Euclidean height at most $T.$ We shall relate this quantity to the
number $N_{\sL}(T)$
of integer Lorentz quadruples of height at most $T$, where
{\em Lorentz quadruples} are those quadruples
that satisfy the Lorentz equation
\beql{302}
-W^2 + X^2 + Y^2 + Z^2  = 0.
\eeq
These are the zero vectors of the {\em Lorentz quadratic form}
\beql{302a}
Q_{\sL}(W, X, Y, Z) = -W^2 + X^2 + Y^2 + Z^2,
\eeq
whose matrix representation is
$$ 
Q_{\sL} =  \left[ \begin{array}{cccc}
          -1 &   0 &   0 &   0 \\
           ~0 &   1 &   0 &   0 \\
           ~0 &   0 &   1 &   0 \\
           ~0 &   0 &   0 &   1 \end{array} \right].
$$
Similarly we shall relate the
number of primitive integer Descartes quadruples, denoted $N_{\sD}^*(T)$,
to the number of primitive integer Lorentz quadruples
of height at most $T$, denoted $N_{\sL}^*(T)$.
We show that there is  a one-to-one height preserving correspondence 
between
integer Descartes quadruples and integer Lorentzian quadruples.
Introduce the matrix $J_0$ defined by
\beql{303}
J_0 = \frac{1}{2} \left[ \begin{array}{cccc}
          1 &  ~1 &  ~1 &  ~1 \\
          1 &  ~1 & -1 & -1 \\
          1 & -1 &  ~1 & -1 \\
          1 & -1 &  -1 &  ~1\end{array} \right]
\eeq
and note that ${J_0}^2 = I.$ The Descartes and Lorentz forms are
related by 
\beql{206a}
 Q_{\sD}=  2 J_0^T Q_{\sL} J_0,
\eeq
which leads to a relation between their zero vectors.

\begin{lemma}~\label{le31}
The mapping $(W,~ X,~  Y,~  Z)^T = J_0 ( w,~  x,~  y, ~ z )^T$ gives a
 bijection
from the set $ ( w,~ x,~ y,~  z )$ of real Descartes quadruples to that of 
real Lorentz quadruples $ (W,~ X,~ Y,~  Z)$ which preserves height. 
It restricts to a
bijection from the set of integer Descartes quadruples to integer Lorentz
quadruples, so that $N_{\sD}(T)= N_{\sL}(T),$ for all $ T > 0,$ and
from primitive integer Descartes quadruples to primitive integer Lorentz
quadruples, so that $N_{\sD}^*(T)= N_{\sL}^*(T),$ for all $ T > 0.$
\end{lemma}

\pf An easy calculation shows that the mapping takes real solutions
of one equation to solutions of the other and that the inverse mapping is
 $( w,~ x,~ y,~ z )^T= J_0 (W,~ X,~ Y,~ Z)^T$, so that it is a bijection.
 The mapping takes integer Descartes quadruples
to integer Lorentz quadruples  because any integer solution to the Descartes
equation satisfies $ w + x + y + z \equiv 0~ (mod~ 2)$. This also holds 
in the reverse direction because integer solutions to the Lorentz 
form also satisfy $W + X + Y + Z \equiv 0~ (mod~ 2)$, as
follows by reducing \eqn{302} $(mod~ 2).$ It is easy to check that
primitive integer Descartes quadruples correspond to primitive 
integer Lorentz quadruples.         $\Box$

Counting the number of integer Descartes quadruples of height 
below a given bound 
$T$ is the same as  
counting integer Lorentz quadruples.
This is a special case of the classical problem 
of estimating the number of representations of a fixed
integer by a fixed diagonal quadratic form, on which there is
an enormous literature. 
For example Ratcliffe and Tschantz~\cite{RT97}
give asymptotic estimates with good error terms
for the number of solutions for the equation $X^2 + Y^2 + Z^2 - W^2 = k,$
of Euclidean height below a given bound, for
all $k \neq 0.$ (They treat  Lorentzian forms in $n$ variables.)
Rather surprisingly the case $k = 0$  seems not to
be treated in the published literature, so 
we derive an asymptotic formula with error term below.
The main term in this asymptotic formula 
was found in 1993 by W. Duke~\cite{Du93}(unpublished)
in the course of establishing an equidistribution result for
its solutions.

\begin{theorem}~\label{th31}
The number of integer Descartes quadruples $N_{\sD}(T)$ of Euclidean height
at most $T$ satisfies $N_{\sD}(T)= N_{\sL}(T)$, and
\beql{304}
  N_{\sL}(T) = \frac {\pi^2 }{4 L(2, \chi_{-4}) }~ T^2 + O (T (log T)^2),
\eeq
as $T \to \infty$, in which 
$$L(2, \chi_{-4}) = \sum_{n = 0}^\infty \frac {(-1)^n}{(2n+1)^2} 
\approx  0.9159.$$
The number $N_{\sD}^*(T)$
of primitive integer Apollonian quadruples of Euclidean height
less than $T$ satisfies $N_{\sD}^*(T)= N_{\sL}^*(T)$ and 
\beql{305}
     N_{\sL}^*(T) = \frac{3} {2 L(2, \chi_{-4})}
~T^2 + O (T(log T)^2),
\eeq
as $T \to \infty.$
\end{theorem}

\pf By Lemma~\ref{le31} it suffices to estimate $N_{\sL}(T).$
Let $r_3 (m)$ denote the number of integer representations
of $m$ as a sum of three squares, allowing positive, negative and
zero integers.  Rewriting the Lorentz
equation as $X^2 + Y^2 + Z^2 = W^2$ we obtain for integer $T$ that
\beql{306}
N_{\sL}( \sqrt{2}T) = 1 + 2 \sum_{m = 1}^T r_3 ( m^2 ) ,
\eeq
since there are two choices for $W$ whenever $W \neq 0.$
A general form for $r_3(m)$ was obtained in 1801 by 
Gauss~\cite[Articles 291-292]{Ga01}, while 
in the special case
$r_3(m^2)$ a simpler form holds, given in 1906  by
Hurwitz~\cite{Hur06}.
This is reformulated in Sandham~\cite[p. 231]{Sa53}
in the form: if $m = \prod_{p} p^{e_p(m)}$, and $p$ runs
over the primes and $m_{odd} = m 2^{-e_2(m)},$ then
\begin{eqnarray}~\label{306a}
r_3(m^2) & =  & 6m_{odd} \prod_{p \equiv 3~(mod~4)}(1 + \frac {2}{p} + .... + 
\frac{2}{p^{e_p(m)}})  \nonumber \\
& = &
6\prod_{ p \equiv 1 ~(mod ~2)} \left( \frac{ p^{e_p(m) + 1} - 1 - 
(\frac {-4}{p})(p^{e_p(m)} - 1)}{p - 1} \right).
\end{eqnarray}
Sandham observes that this formula is equivalent to
\beql{306b}
F(s) := \sum_{m = 1}^\infty \frac {r_3 (m^2)}{m^{s}} = 
 6 (1 - 2^{1-s}) \frac {\zeta(s) \zeta(s - 1)}{L (s, \chi_{-4})}
\eeq
where 
\beql{306c}
L (s, \chi_{-4}) := \sum_{m = 1}^\infty (\frac {-4}{m}) m^{-s} 
= \sum_{n = 0}^\infty \frac {(-1)^n}{(2n + 1)^s}.
\eeq
The right hand side of \eqn{306b} is a meromorphic function in the $s$-plane,
 which has a simple pole at $s = 2$ with residue 
$$c_1 = \frac {3 \zeta(2)}{L(2, \chi_{-4})} = 
\frac { \pi^2}{2 L(2, \chi_{-4}) },$$ 
and has no other poles for $\Re{s} > 1.$ 
One could now proceed by a  standard contour
integral approach~\footnote{Integrate $F(s)$ against $\frac{T^s}{s} ds$ on
a vertical contour $\Re{s} = c$ for some $c$ with $1 < c < 2$,
cf. Davenport~\cite[Chapter 17]{Da00}.}
to obtain $N_{\sL}(\sqrt{2}T)$ equals 
 $c_1 T^2$  plus a lower-order error
term, so that $N_{\sL}(T) =  \frac{1}{2} c_1 T^2 + o(T^2).$
However we will derive \eqn{304} 
directly from
the exact formula \eqn{306a}. We have, 
\begin{eqnarray*}
N_{\sL}(\sqrt{2}T) & =  & 1 +  2\sum_{1 \leq j \leq log_2 T} 
\sum_{n = 1}^{\lfloor \frac{T}{2^j} + \frac{1}{2}\rfloor} r_3 ( (2n - 1)^2) \\
& = &
1 + 12 \sum_{1 \leq j \leq log_2 T} 
\sum_{n = 1}^{\lfloor\frac{T}{2^j}+ \frac{1}{2} \rfloor} (2n - 1) 
\prod_{{p \equiv 3~ (mod~ 4)}\atop{p | 2n - 1}} 
(1 + \frac {2}{p} + .... + \frac{2}{p^{e_p(2n - 1)}}).
\end{eqnarray*}
Expanding the products above and using
 $\sum_{n = 1}^U (2n - 1) = U^2 + O (U)$,
one obtains
\begin{eqnarray}~\label{308}
N_{\sL}(\sqrt{2}T) &= &
1 + 12 \sum_{k \ge 0} \sum_{P_k} \frac{2^k}{P_k} \left(\sum_{j \ge 1}
\sum_{m=1}^{\lfloor \frac {T}{2^j P_k}+ \frac{1}{2}\rfloor} (2m-1) P_k \right)
\nonumber \\
&= & \frac {12}{4} \left(
\sum_{k \geq 0} 2^k \sum_{{{j\ge 0, P_k}\atop {P_k < T/2^j}}}
 \frac{1}{2^{2j} P_k^{2}} \right) T^2 + 
O \left( (\sum_{k \geq 0} 2^k \sum_{P_k \le T} \frac{1}{P_k})T \right),
\end{eqnarray}
in which $P_k = p_1^{e_1}...p_k^{e_k}$
with all $p_i \equiv ~3~(mod~4)$ and all $e_i \geq 1.$ 
If the condition $P_k < T/ 2^j$ were dropped in the first
sum in parentheses above, then
it would sum to $\frac {\zeta(2)}{L(2, \chi_{-4})},$
as one sees by examining the associated Euler product, which is
$$ 
\frac {\zeta(s)}{L(s, \chi_{-4})}= 
(1 - 2^{-s})^{-1} \prod_{p \equiv 3~(mod~4)} \frac {1 + p^{-s}}{1 - p^{-s}}
= \sum_{m=1}^\infty \frac{a(m)}{m^s}  ,$$
evaluated at $s = 2$, since
$\frac {1 + p^{-s}}{1 - p^{-s}} = 1 + 2 p^{-s} + 2 p^{-2s} + ...$
 The error introduced by truncating this Dirichlet series 
at $s=2$ at $n \le  T$ 
is  bounded by 
$$
\sum_{m=T}^\infty  \frac{a(m)}{ m^2} \le
\sum_{m=T}^\infty \frac {2^{\nu(m)}}{m^2} \le 
\sum_{m=T}^\infty \frac {d(m)}{m^2}  = O(\frac{\log T}{T}),
$$ 
in which $\nu(n)$ counts the number of distinct prime divisors of $m$
(without multiplicity), $d(n)$ counts the number of divisors on $n$,
and the last estimate uses 
\footnote{ Sum in blocks $2^j T \le n \le 2^{j+1}T$, viewing
denominators as nearly constant and averaging over numerators.}
the fact 
(\cite[Theorem 319]{HW}) that
the average order of $d(n)$ is $\log n$.
The remainder term in \eqn{308} is bounded by $O (T (\log T)^2)$, because
$$
\sum_{k \geq 0} 2^k \sum_{P_k \le T} \frac{1}{P_k} \le 
\sum_{m=1}^T \frac{ 2^{\nu(m)}}{m} \le \sum_{m=1}^T \frac{d(m)}{m}
\le ( \sum_{m=1}^T \frac{1}{m})^2     = O ((\log T)^2 ).
$$
Since $\zeta (2) = \frac {\pi^2}{6},$ these estimates
combine to give
$$N_{\sL}(\sqrt{2}T) = \frac {\pi^2}{2L(2, \chi_{-4})} T^2 + 
O(T~(log T)^2).$$
Rescaling $T$ yields \eqn{304}.

To handle primitive Lorentz quadruples, we use the
function $r_3^* (m)$ which counts the number of primitive integer
representations of $m$ as a sum of three squares, using positive,
negative and zero integers. Then
$r_3(m^2) = \sum_{d~|~m} r_3^* (d^2)$, which
yields
$$
F^{\ast}(s) := \sum_{m = 1}^\infty \frac {r_3^{\ast} (m^2)}{m^{s}} = 
 \frac{F(s)}{\zeta (s)} =
 6 (1 - 2^{1-s}) \frac {\zeta(s - 1)}{L (s, \chi_{-4})}.
$$
The function  $F^{\ast}(s)$ is analytic in $\Re(s) >1$ 
except for a simple  pole at $s=2$ with residue
$\frac{6}{\pi^2} c_1$, which gives the
constant in the main term of \eqn{305}.
 To obtain the
error estimate \eqn{305} one can use \eqn{304} and M\"obius inversion,
with $r_3^*(m^2) = \sum_{d~|~m} \mu(d) r_3( (\frac{m}{d})^2)$.
We omit details.
$\Box$
 
\paragraph{Remarks.} 
(1) Various Dirichlet series associated to zero solutions of
indefinite quadratic forms have meromorphic continuations to $\cc$,
cf. Andrianov~\cite{And89}. These can used to obtain asymptotics
for the number of solutions satisfying  various side conditions.

(2) The real solutions of the homogeneous equation 
$Q_{\sL}(w, x,y, z)= -w^2 + x^2 +y^2 +z^2 = 0$  form the {\em light cone} in
special relativity.

%--------------------------------------------------------------------
%
% Section 3
%
%--------------------------------------------------------------------
%

\section{Reduction Theory and Root Quadruples}
\setcounter{equation}{0}
In this section we describe, given an  Apollonian circle packing
and a  Descartes quadruple in it, 
a reduction procedure which, if it halts, identifies
within it a unique Descartes quadruple $(a, b, c, d)$ 
which is ``minimal''. 
This quadruple is called the {\em root quadruple} of the packing. 
This procedure always halts for integral Apollonian packings.

\begin{defi}
{\em
The {\em Apollonian group} $\sA$ is the group generated by the four
integer $4 \times 4$ matrices
\beqs
S_1=\left[ \begin{array}{cccc}
          -1 & 2 & 2 & 2 \\
          0 & 1 & 0 & 0 \\
          0 & 0 & 1 & 0 \\
          0 & 0 & 0 & 1\end{array} \right]
\qquad
S_2=\left[ \begin{array}{cccc}
          1 & 0 & 0 & 0 \\
          2 & -1 & 2 & 2 \\
          0 & 0 & 1 & 0 \\
          0 & 0 & 0 & 1\end{array} \right]
\eeqs
\beqs
\qquad
S_3=\left[ \begin{array}{cccc}
          1 & 0 & 0 & 0 \\
          0 & 1 & 0 & 0 \\
          2 & 2 & -1 & 2 \\
          0 & 0 & 0 & 1\end{array} \right]
\qquad
S_4=\left[ \begin{array}{cccc}
          1 & 0 & 0 & 0 \\
          0 & 1 & 0 & 0 \\
          0 & 0 & 1 & 0 \\
          2 & 2 & 2 & -1\end{array} \right]
\eeqs
}
\end{defi}

As mentioned earlier, the Apollonian group was introduced in
the 1967 paper of Hirst\cite{Hirst67}, and 
was later used in  S\"{o}derberg \cite{So92} and Aharonov and
Stephenson~\cite{AS97} in studying Apollonian packings.

We view real Descartes quadruples $\bv= (a,~b,~c,~d)^T$ as column
vectors, and the Apollonian group $\sA$ acts by matrix 
multiplication, sending $\bv$ to $M\bv$, for $M \in \sA$.
 The action takes Descartes quadruples to Descartes
quadruples, because $\sA \subset Aut_{\zz}(Q_{\sD})$,
the set of real automorphs of 
the Descartes integral quadratic form $Q_{\sD}$ given in
  \eqn{300a}. That is, each such $M$ satisfies
$$ 
M^T Q_{\sD} M = Q_{\sD},  \qquad\mbox{for~all}~ M \in \sA,
$$
a relation which it suffices to check on the four generators
$S_i \in \sA.$

The elements $S_j$ have a geometric meaning as corresponding
to inversion in one of the four circles of a Descartes quadruple
to give a new quadruple in the same circle packing, as explained
in \cite[Section 2]{GLMWY21}. That paper showed
that this group with the given
generators is a Coxeter group whose only relations are
$S_1^2=S_2^2=S_3^2=S_4^2=I$.

The reduction procedure attempts to reduce the size of the elements
in a Descartes quadruple by applying one of the generators $S_j$
to take the quadruple $\bv = (a,~b,~c,~d)$ viewed as a column vector
to the new quadruple $S_j \bv,$ until further decrease is not
possible. Each packing has a well-defined {\em sign}, which is
the sign of $L(\bv):=a+b+c+d$ of any Descartes quadruple in the 
packing. (The
sign is independent of the quadruple chosen in the packing; see 
Lemma~\ref{le31a}(ii) below.)
We describe the reduction procedure
for  packings with positive sign, those with
$L(\bv) := a + b +c +d > 0$; the reduction procedure for negative
sign packings is obtained by conjugating by  the inversion
$(a, b,c,d) \to (-a,-b,-c,-d)$, which takes negative sign packings
to positive sign packings.
Thus we  suppose $L(\bv) := a + b +c +d > 0$ and 
for simplicity consider the case where the quadruple is ordered
$a \le b \le c \le d$. We consider which $S_j \bv$ can decrease the
sum $L(\bv) = a + b + c + d$, and show below this can never occur
using $S_1, S_2$ or $S_3$. Note that 
$S_4 (a,~b,~c,~d)^T = (a,~b,~c,~d')^T$ where
$d' = 2(a + b + c) - d.$
 
\begin{lemma}\label{le31a}
Suppose that $\bv = (a,~b,~c,~d)^T$ is
a real Descartes quadruple.
Let $\bv$ 
have its elements ordered
$a \leq b \leq c \leq d$, and set $d' = 2(a + b+ c) - d$.

(i) If  $L(\bv)= a+b+c+d>0$, 
then $a+b \geq 0$, with equality holding only if $a=b=0$ and $c=d$. 
As a consequence, we always have $b\geq 0$. 

(ii) If  $L(\bv)= a+ b+ c+ d > 0$, then 
each  $L(S_j\bv) >0.$ In particular      $L(S_4\bv) = a+ b+ c+ d'>0.$

(iii) If $a \geq 0 $, so that $L(\bv)= a + b + c + d \geq 0$, 
then $d' \leq c \leq d. $ 
If $d' < c$ then the matrix  $S_4$ that changes $d$ to $d'$
strictly decreases the sum $L(\bv):= a+b+c+d,$
and it is the only generator of $\sA$ that does so. If $d' = c$
then necessarily $c = d = d'$ and no generator $S_i$ of $\sA$
decreases $L(\bv)$.

\end{lemma}
\pf
(i) If $a \geq  0 $ then we are done, so assume $a < 0$. 
Suppose first that $0 \leq b \leq c\leq d$. 
Let $x=-(a+b)$, $y=-ab$, so that  $y\geq 0$.
In terms of these variables, the  Descartes equation becomes
$$
 2(x^2+2y+c^2+d^2)=(-x + c+d)^2, 
$$
which simplifies to 
\begin{equation}~\label{three}
 (c-d)^2+x^2+4y+2(c+d)x = 0.
\end{equation}
This equation cannot hold if $x>0$, because all 
terms on the left are nonnegative, and some are positive, hence
$a+b \ge 0$. If 
$x=0$, then \eqref{three} implies 
$y=0$ and $c=d$, whence $a=b=0$ and $c=d$. 

The remaining case is  $a \leq b < 0.$ We first observe 
that in any Descartes quadruple, at least three terms have
the same sign. Indeed the Descartes equation
can be rearranged as
$$ 
(a-b)^2 + (c-d)^2= (a+b)(c+d),
$$
and if 
$ a \le b < 0 < c \le d$
then the left side is nonnegative while the right side is
strictly negative, a contradiction.
Thus we must have $a \leq b \leq c \leq 0 < d$,
since $a + b+c+d >0$.
Now the Descartes quadruple
$(-d, -c, -b, -a)$ has
$ -d < 0 \le -c $, so the  preceding case 
applies to show  that $-d-c \geq 0$. 
We conclude that  $a+b+c+d \leq 0$, contradicting the fact that $a+b+c+d >0$. 
This completes the proof of (i).

(ii) The ordering of entries gives $L(S_1\bv) \ge L(S_2\bv) \ge L(S_3\bv) \ge
L(S_4\bv)$, so it suffices to prove $L(S_4\bv) > 0$.
The Descartes equation implies that
\[
d, d'=a+b+c\pm 2q_{abc}, \qquad \text{where} \qquad q_{abc}=\sqrt{ab+bc+ca}.
\]
We have
$a+b+c+d'=2(a+b+c)-2\sqrt{ab+bc+ac} > 0$ because $a+b+c \geq 0$ 
(using (i)) and
\[
(a+b+c)^2-(ab+bc+ac)=\frac{1}{2} ((a+b)^2+(b+c)^2+(a+c)^2) > 0. 
\]

(iii) The Descartes equation \eqref{descartes} gives 
\[
d'=a+b+c-2\sqrt{ab+bc+ac}.
\]
Thus
\[
d'-c=a+b-2\sqrt{ab+bc+ac} \leq a+b-\sqrt{4(a+b)c} \leq a+b-\sqrt{(a+b)^2}=0.
\]

If $d' < c \le d$ then the sum $ L(\bv')= a+b+c+d' < L(\bv),$ so
the sum decreases. If $S_i$ changes
$c$ to $c'$, then $c' = 2(a + b + d) - c \geq 2(a+b+c) - c \geq c$
because $a + b \geq 0$ by (i), so the sum $L(\bv)$ does not decrease in
this case. Similarly the sum does not decrease if $S_i$ changes
$b$ to $b'$ or $a$ to $a'.$
In the case of equality $d' = c$, one easily checks that  $c = d = d'$ ,
which forces $a = b = 0,$ and no $S_i$ decreases the
sum $L(\bv).$ 
$\Box$

\begin{defi}
{\em A Descartes quadruple $\bv=(a, b, c, d)$ with $L(\bv)= a+b+c+d>0$ is a 
{\em root quadruple } if 
$a \leq 0 \leq b\leq c \leq d$ and $a+b+c \geq d$.}
\end{defi} 
 
Note that the  last inequality above 
is  equivalent to the condition  $d' = 2(a+b+c) - d \geq d$. \\

\noindent {\bf Reduction algorithm.} \\
{\em Input: A real Descartes quadruple $(a,~b,~c,~d)$ with
$a + b+ c+ d > 0$.

(1) Test in order $1 \le i \le 4$ whether some $S_i$
decreases the sum $a+b+c+d.$  It so, apply it to produce a
new quadruple and continue. 

(2) If no $S_i$ decreases the sum, order the elements of the
quadruple in increasing order and halt.} \\

The reduction algorithm takes real quadruples as input, and is
not always guaranteed to halt. The following theorem shows
that when the algorithm is given
an integer Descartes quadruple as input, it always halts, and 
outputs a root quadruple. In the algorithm, 
the element $S_i$ that decreases the sum necessarily
decreases the largest element in the quadruple, leaving the
other three elements unchanged. The proof below 
establishes that in all cases where
a reduction is possible, the largest element of
the quadruple is unique, so that the choice of $S_i$ in the
reduction step is unique. There do exist quadruples with
a tie in the largest element, such as $(0,~0,~1,~1),$ but the
vector $(a,~b,~c,~d)$ cannot then be further reduced. 

\begin{theorem}~\label{th41b}
{\rm (1)} If the reduction algorithm ever encounters some element $a<0$, 
then it will halt at a root quadruple in finitely many more steps.

{\rm (2)} If $a, b, c, d$ are  integers, then the reduction 
algorithm will halt  at a root quadruple in finitely many steps.
 
{\rm (3)} A root quadruple is unique if it exists.
However an Apollonian circle packing 
may contain more than one Descartes configuration yielding this quadruple. 
\end{theorem}

\pf

(1) Geometrically a Descartes quadruple with $a<0$ describes a circle of
 of radius $1/a$ enclosing three
mutually tangent circles of radii $1/b, 1/c, 1/d$. All circles
in the packing lie inside this bounding circle of radius  $1/a$.
Each non-halting reduction produces
a new circle of radius $1/d' > 1/d$, which covers an area of 
$\pi/{d'}^2$, and this is at least $\pi/d^2.$ 
Since there is a total area of $\pi /a^2$ which can be covered, 
and all circles except the one with radius $1/a$ have disjoint interiors, 
this process must halt in at most $\displaystyle
\Big\lf (\frac{d}{a})^2\Big\rf$ steps.   

(2)
Let $q_{a b c}=\sqrt{ab+bc+ac}=(a+b+c-d)/2 \in \nn$. 
After each reduction, the sum
$a+b+c+d$ decreases by $4q_{abc}$. 
By Lemma \ref{le31a}, the sum $a+b+c+d$ is bounded below by
$0$. Therefore this process halts after finitely many steps.

(3) 
If $(a, b, c, d)$ is a root quadruple of an
Apollonian packing,  
then the numbers $a,b, c, d$ are the curvatures of the largest
circles contained in this packing, hence they are unique. 
%[ We haven't yet defined double semi-bounded, have we?]
On the other hand, the Apollonian packing may contain more than one
copy of this quadruple, for example, $(-1, 2, 2, 3)$ appears twice in
the packing shown in Figure \ref{std}, and $(0, 0, 1,1 )$ appears
infinitely many times  in the packing in 
Figure~\ref{vertical}.
$\Box$

Root quadruples lead to a partition of the set $Q( \zz )$ of all 
integer Descartes quadruples.
This set partitions into $ Q( \zz )^+ \cup \{(0,0,0,0)\} \cup  Q( \zz )^-,$
where
\beql{CI32}
Q( \zz )^+ = \{ (a,b,c,d) \in Q( \zz ) : a +b+c+d > 0 \}
\eeq
and $Q( \zz )^- = -Q(\zz)^+$.
Next we  have the partition
\beql{CI33}
Q(\zz )^+ = \bigcup_{k=1}^\infty kQ(\zz )_{prim}^+ ~,
\eeq
where $Q( \zz )_{prim}^+$ enumerates all primitive integer Descartes 
quadruples in $Q(\zz )^+.$
These latter are exactly the Descartes quadruples occurring in all primitive 
integer Apollonian packings,
so we may further partition $Q(\zz )_{prim}^+$ into 
a union of the sets $Q( \sP_{\sD} )$, where  $Q( \sP_{\sD} )$    denotes
the set of all Descartes quadruples in the circle packing
$\sP_{\sD}$ having primitive root
quadruple $\sD$, i.e.
\beql{CI34}
Q(\zz )_{prim}^+ = \bigcup_{\begin{array}{c}
\mbox{primitive root} \\
\mbox{quadruple $\sD$}
\end{array}
}
Q( \sP_{\sD} ).
\eeq
We study the distribution of root quadruples in \S4 and the set of
integers in a given packing $\sP_{\sD}$ in \S5 and \S6.

By definition the  Apollonian group labels all
the (unordered) Descartes quadruples
in a fixed Apollonian packing. We now show that it has the
additional property that 
for a given integral Apollonian
packing, the integer curvatures
of all circles not in the root quadruple lie
in one-to-one correspondence with the non-identity elements of
the Apollonian group.

\begin{theorem}~\label{th42b}
Let $\sP_{\bv}$ be the integer Apollonian circle packing with root quadruple
$\bv = (a,~b,~c,~d)^T$.
%  and suppose that $a < 0$.
 Then the set of integer curvatures occurring
in $\sP$, counted with multiplicity, consists of the four elements of
$\bv$ plus the largest elements of each vector $M\bv$, where $M$ runs over all
elements of the Apollonian group $\sA.$
\end{theorem}

\paragraph{Proof.}
%The condition $a <0$ excludes the packing $(0,0,1,1)$.
Let
$M= S_{i_n} \cdots S_{i_1}$ be a reduced word in the generators of $\sA$,
that is $S_{i_k} \neq S_{i_{k+1}}$ for $1 \le k < n.$ The main
point of the proof is that if
$\bw^{(n)} = S_{i_n} \cdots S_{i_1} \bv$,
then $\bw^{(n)}$ is obtained from $\bw^{(n-1)}$ by changing one
entry, and the new entry inserted is always the largest entry in the 
new vector. (It may be tied for largest value.)
We prove this by induction on $n$.
In the base case $n=1$,
there are four possible vectors $S_i \bv$, whose inserted entries are
$a'= 2(b+c+d)-a$,
$b' = 2(a+c+d) -b$, $c' = 2(a+b+d) -c$, and $d' = 2(a+b+c) -d$,
respectively.
Since $a \le b \le c \le d$ we have $d' \le c ' \le b' \le a'$,
and since $\bv$ is a root quadruple, we have $d' \geq d$,
as asserted.

For the induction step, where $n \ge 2$, there are only three choices for 
$S_{i_n}$ since $S_{i_n} \neq S_{i_{n-1}}$.
If the elements of $\bw^{(n-1)}$ are labelled in increasing order as 
$w_1^{(n-1)} \le w_2^{(n-1)} \le w_3^{(n-1)} \le w_4^{(n-1)}$, 
then we may choose the labels (in case of a tie for the largest
element) so that $w_4^{(n-1)}$ was
produced at step $n-1$, by the induction hypothesis.
Thus exchanging $w_4^{(n-1)}$ is forbidden at step $n$, hence if
$w_4^{(n)}$ denotes the new value produced at the next step, then
\begin{eqnarray}\label{DI32}
w_4^{(n)} & \ge & 2(w_1^{(n-1)} + w_2^{(n-1)} + w_4^{(n-1)} - w_3^{(n-1)} \nonumber \\
& \ge & 2w_1^{(n-1)} + 2w_2^{(n-1)} + w_4^{(n-1)} > w_4^{(n-1)},
\end{eqnarray}
because $w_1^{(n-1)} + w_2^{(n-1)} > 0$ by Lemma \ref{le31a}(i).
This completes the induction step.

The inversion operation produces one new circle in the packing,
namely the new value added in the Descartes quadruple, and \eqn{DI32}
shows that its curvature is
$||M \bv||_\infty$, where 
$||\bw ||_\infty $ denotes the supremum norm of the vector $\bw$.

Every circle in the packing is produced in this procedure, 
by definition of the Apollonian group.
That all non-identity words $M \in \sA$ label distinct circles is 
clear geometrically from
the tree structure of the packing.~~~$\Box$

%--------------------------------------------------------------------
%
% Section 4
%
%--------------------------------------------------------------------
%

\section{Distribution of Primitive Integer Root Quadruples}
\setcounter{equation}{0}
In this section we count 
integer Apollonian circle packings
in terms of the size of their root quadruples. Recall that
a Descartes quadruple $(a,~b,~c,~d)$ is a {\em root quadruple}
if $a \leq 0 \leq b \leq c \leq d$ and $ d' = 2(a + b + c) - d \geq d > 0.$
It suffices to study
primitive packings, i.e. ones whose integer quadruples are
relatively prime. We begin with a Diophantine characterization
of root quadruples.

\begin{theorem}\label{th51}
Given a solution $(a,b,c,d) \in \zz^4$ 
to the Descartes  equation 
$$(a+b+c+d)^2 = 2 (a^2 + b^2 + c^2 + d^2 ),$$
define $(x,d_1, d_2, m)$ by
\beql{N51}
\left[\begin{array}{c}
a \\ b \\ c \\ d
\end{array} \right] =
\left[
\begin{array}{rccr}
1 & 0 & 0 & 0 \\
-1 & 1 & 0 & 0 \\
-1 & 0 & 1 & 0 \\
-1 & 1 & 1 & -2
\end{array}
\right]
\left[\begin{array}{c}
x \\ d_1 \\ d_2 \\ m
\end{array} \right] =
\left[ \begin{array}{c}
x \\ d_1 -x \\ d_2 -x \\ -2m+d_1 + d_2 -x
\end{array}
\right] \,.
\eeq
Then $(x,d_1,d_2, m ) \in \zz^4$ satisfies
\beql{N51a}
x^2 + m^2 = d_1 d_2 ~.
\eeq
Conversely, any solution $(x, d_1, d_2, m) \in \zz^4$ to this
equation yields an integer solution
to the Descartes equation as above.
In addition:
\begin{itemize}
\item[(i)]
The solution $(a,b,c,d)$ is primitive if and only if 
$\gcd (x, d_1, d_2 ) =1$.
\item[(ii)]
The solution $(a,b,c,d)$ with
$a+b+c+d >0$ is a root quadruple  if
and only if \\
$x < 0 \le 2m \le d_1 \le d_2$.
\end{itemize}
\end{theorem}

\paragraph{Proof.}
The first part of the theorem requires, to have $m \in \zz$, that
$ a+b+c+d \equiv 0~(\bmod~2)$.
This follows from the Descartes equation by reduction $(\bmod~2)$.

For (i), note that $\gcd (x,d_1, d_2) = \gcd (a,b,c) = \gcd (a,b,c,d)$.

For (ii), the condition $a < 0 \le b \le c \le d$ implies successively $x< 0$,
$d_1 \le d_2$, $d_1 - 2x = b-a \ge 0$, and $-2m +d_1 = d-c \ge 0$.
Finally the root quadruple
condition $d' = 2(a+b+c)-d \ge d \ge 0$ gives $d' = 2m \ge 0$.
Thus $x < 0 \le 2m \le d_1 \le d_2$.
The converse implication follows similarly.~~~$\Box$

We now study primitive integer root quadruples with $a= -n$, 
for $n \in \zz_{\ge 0}$.
Theorem \ref{th51} shows that they are in one-to-one correspondence with the 
integer solutions $(m, d_1, d_2 )$ to
\beql{eq14}
n^2 + m^2 = d_1 d_2,
\eeq
\beql{eq15}
0 \le 2m \le d_1 \le d_2 \quad\mbox{and}\quad
\gcd (n, d_1, d_2 ) =1 ~.
\eeq
For each of $n=0,1,2$, there is  one primitive root quadruple 
with $a = -n$, 
namely, $(0,0,1,1)$, $(-1,2,2,3)$, $(-2,3,6,7)$, respectively.
For $n=3$, there are two, $(-3,4,12,13)$ and $(-3,5,8,8)$.
As an example of a nonsymmetric integral Apollonian circle
packing, Figure 3 pictures the packing $(-6,~11,~14,~15).$

\begin{figure}[htbp]\label{nonsymmetric}
\centerline{\epsfxsize=4.0in \epsfbox{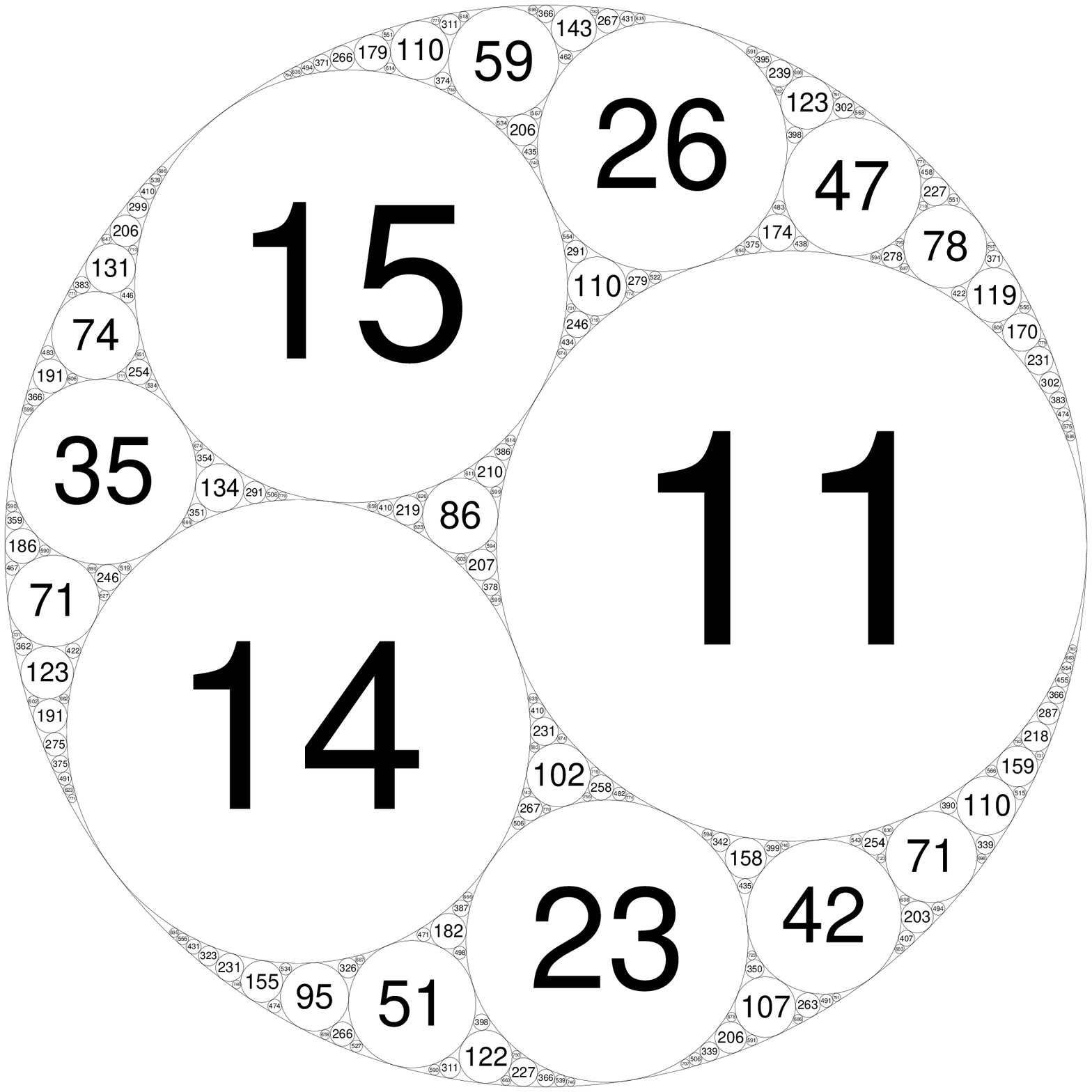}}
\caption{The nonsymmetric packing $(-6, ~11,~14,~15).$}
\end{figure}

Let $N_{root}^*(-n)$ denote the number of primitive root quadruples
with negative element $-n$.
One has $N_{root}^*(-n) \ge 1$ for all $n \ge 0$, since
$(x, ~d_1, ~d_2, ~m) = (-n, 1, n^2, 0)$ in Theorem~\ref{th51} produces
the primitive root quadruple 
$(a,~b,~c, ~d)= (-n, n+1, n(n+ 1), n(n+1) + 1 )$. 
Table 1 and Table 2 present selected
values of $N_{root}^*(-n)$.

\begin{table}[htpb]
\begin{center}
\begin{tabular}{|c|c||c|c||c|c||c|c||c|c||} \hline
$n$ & $N(-n)$ & $n$ & $N(-n)$ & $n$ & $N(-n)$ & $n$ & $N(-n)$ &$n$ & $N(-n)$ \\
\hline \hline
1 & 1 &  11 & 4 &  21 & 10  &  31 & 9 &  41 & 11 \\ \hline
2 & 1 &  12 & 6 &  22 &  7  &  32 & 9 &  42 & 18 \\ \hline
3 & 2 &  13 & 4 &  23 &  7  &  33 & 14&  43 & 12 \\ \hline
4 & 2 &  14 & 5 &  24 & 10  &  34 & 9 &  44 & 14 \\ \hline
5 & 2 &  15 & 6 &  25 &  6  &  35 & 10&  45 & 14 \\ \hline
6 & 3 &  16 & 5 &  26 &  7  &  36 & 14&  46 & 13 \\ \hline
7 & 3 &  17 & 5 &  27 & 10  &  37 & 10&  47 & 13 \\ \hline
8 & 3 &  18 & 7 &  28 & 10  &  38 & 11&  48 & 18 \\ \hline
9 & 4 &  19 & 6 &  29 &  8  &  39 & 14&  49 & 15 \\ \hline
10& 3 &  20 & 6 &  30 & 10  &  40 & 10&  50 & 11 \\ \hline
\end{tabular}
\end{center}
\caption{$N_{root}^*(-n)$ for small $n$}
\label{tab51}
\end{table}

\begin{table}[htpb]
\begin{center}
\begin{tabular}{|c|c||c|c||c|c||c|c||} \hline
$p$ &  $N(-p)$ & $p$ & $N(-p)$ & $p$  & $ N(-p)$ & $p$ & $N(-p)$ \\
\hline\hline
1009 & 253 &  3001  & 751  & 4007 & 1003 &  5011 & 1254 \\
\hline
1013 & 254 &  3011  & 754  & 4013 & 1004 & 10007 & 2503 \\
\hline
2003 & 502 &  4001  & 1001 & 5003 & 1252 & 10009 & 2503 \\
\hline
2011 & 504 &  4003  & 1002 & 5009 & 1253 & 20011 & 5004 \\
\hline
\end{tabular}
\end{center}
\caption{$N_{root}^*(-p)$ for selected prime $p$.}
\label{tab52}
\end{table}

Examination of  numerical data 
led the the third author
and S. Northshield independently to conjecture an exact formula
for $N_{root}^*(-n)$, 
stated as  Theorem~\ref{Nth42} below. Northshield~\cite{No02}
obtained a proof, and we give another here.
We establish this result by showing that
$N_{root}^*(-n)$ can be interpreted as a class number, 
namely   $ N_{root}^*(-n)  = h^{\pm}(-4n^2)$
where for $-4\Delta < 0$ the quantity $h^{\pm}(-4\Delta)$ counts the number
of $GL(2,\zz)$-equivalence classes of primitive integral
binary quadratic forms 
$$[A, B, C]:= AT^2+2BTU +CU^2,$$ 
of discriminant $-4\Delta = 4B^2 - 4AC$; a form is primitive
if $gcd(A, B, C) =1$. By
comparing $GL(2,\zz)$-actions and $SL(2,\zz)$-actions
on positive definite binary quadratic forms one obtains  the relation
\beql{N200b}
h^{\pm}(-4\Delta) = \frac{1}{2}( h(-4\Delta) + a(-4\Delta))
\eeq
where $h(-4\Delta)$ is the usual binary quadratic
form class number (for $SL(2, \zz)$-equivalence) and
$a(-4\Delta)$ counts the $SL(2, \zz)$-equivalence classes containing
an ambiguous reduced binary form of discriminant
$-4\Delta$ (defined below).
The exact formula in Theorem~\ref{Nth42} 
is deduced from  classical formulas for $h(-4n^2)$ and
$a(-4n^2)$.

\begin{theorem}~\label{Nth41}
The Apollonian root quadruples $(-n, x,y,z)$
with $-n < 0 \le x \le y \le z$ are in one-to-one correspondence with
 positive definite integral binary quadratic forms of discriminant
$ -4n^2$ having nonnegative middle
coefficient.  The 
associated reduced binary  quadratic
form $[A, B,C] = A T^2 +2BTU +CU^2$ is given by
$$[A, B, C]:= [-n+x,~ \frac{1}{2}(-n+x+y-z),~ -n+y]. $$
Primitive root quadruples correspond to reduced binary quadratic
forms having a nonnegative middle coefficient.
In particular, the  number  $N_{root}^*(-n)$
 of primitive root quadruples with least element $-n$
satisfiies
$$
N_{root}^*(-n) = h^{\pm}(-4n^2),
$$
where
$ h^{\pm}(-4n^2)$ is the number of $GL(2,\zz)$-equivalence
classes of positive definite primitive binary integral forms
of discriminant $-4n^2$.
\end{theorem} 

\noindent\paragraph{Remark.}
A positive definite form $[A, B, C]$ is {\em reduced} if 
$0 \le |2B| \le A \le C.$ 
Reduced forms with nonnegative middle
coefficients enumerate form clases under $GL(2, \zz)$-action,
rather than the $SL(2, \zz)$-action studied by
Gauss; the $GL(2,\zz)$-action
corresponds to the notion of 
equivalence of (definite) quadratic forms used by Legendre.\\

\noindent{\bf Proof.}
Recall that the  Descartes quadratic form is
$$Q_{\sD}(w,x,y,z) := (w+x+y+z)^2 - 2(w^2 +x^2+y^2+z^2).$$
A Descartes quadruple is any integer solution to 
$Q_{\sD}(w,x,y,z) =0$
and  a Descartes quadruple is 
 {\em primitive} if $gcd(w,x,y,z)=1.$ Any 
integer solution  satisfies
$$
w+x+y+z \equiv ~0~(\bmod~2),
$$
as follows by reduction $(\bmod~2)$.
By definition a Descartes quadruple $(w, x, y, z)$ satisfying
$w+x+y+z >0$ is a {\em root quadruple} if and only if 

(1)  $w= -n \le  0 \le x \le y \le z$,

(2) $ 2(-n +x+y)-z \ge z.$ 

\noindent Condition (1) implies conversely that
$w+x+y+z >0$, and condition (2) is equivalent to
$$ -n + x + y - z \ge 0.$$

The integer  solutions $Q_{\sD}(w,x,y,z)=0$ with $w = -n <0$ are in one-to-one
correspondence with integer representations of $n^2$ by the
the ternary quadratic form 
$$Q_T(X,Y,Z) := XY + XZ + YZ,$$
The correspondence between solutions to $Q_{\sD}(w,x,y,z)=0$
and solutions to $Q_T(X,Y,Z)=n^2$ is given by
$$
(X, Y, Z) := (\frac{1}{2}(w+x+y - z),\frac{1}{2}(w + x -y + z),
\frac{1}{2}(w-x+y+z)).
$$
The congruence condition
$w + x + y + z \equiv 0 ~(\bmod~2)$
on integer Descartes quadruples 
implies that $(X, Y, Z)$  are all integers. 
The map is onto, because an explicit inverse map from
an integer solution  $(X,Y, Z)$    is
$$
(w, x, y, z) := (-n, n+X+Y, n+X+Z, n+Y+Z).
$$
The primitivity condition $gcd(-n, x,y,z)=1$ on
Descartes quadruples is equivalent to  the
primitivity condition $gcd(X, Y, Z)=1$. 
The ``root quadruple conditions''  (1) and (2) above are
equivalent to the inequalities
$$0 \le X \le Y \le Z.$$

For any integer $M$,
the integer solutions of $Q_T (X,Y,Z) =M$ are in one-to-one correspondence
with integer representations of $M$ by the 
determinant ternary quadratic form 
\beql{N200c}
Q_\Delta(A, B, C):= AC- B^2. 
\eeq
A solution $(X,Y, Z)$ gives a solution to 
$Q_\Delta(A, B, C) = M$
under the substitution 
$$
(A,B,C) := (X+Y, X,X+Z).
$$
The inverse map is
$$
(X,Y, Z) :=(B, A-B, C-B).
$$
The ``root quadruple  conditions''  above on $(X, Y, Z)$ are easily checked
to be  equivalent to the inequalities
\beql{N201}
 0 \le 2B \le A \le C.
\eeq
In this correspondence.
$B$ is necessarily  an integer, so  $2B$ is an even integer.
The primitivity condition $gcd(X, Y, Z)=1$ transforms to the
primitivity condition $gcd(A, B, C)=1.$

The conditions \eqn{N201} give a complete set of
equivalence classes of primitive positive definite integral binary forms 
 of fixed discriminant  $D= -4M$   under the action of $GL(2,\zz).$
To show this, note that the conditions for the positive  definite 
binary quadratic form
$[A, B, C] := A T^2 + 2B TU + C U^2$ of discriminant 
$D = 4B^2 - 4AC = -4M$
to be a reduced form in the $SL(2,\zz)$ sense are
\beql{N202}
 0 \le |2B| \le A \le C.
\eeq
It is known that all positive
definite integral forms are $SL(2,\zz)$-equivalent 
to a reduced form, and that all reduced forms are 
$SL(2, \zz)$-inequivalent, except for 
$[A, B, A] \approx [A,-B, A]$ and $[2A, A, C] \approx [2A, -A, C]$,
where $\approx$ denotes $SL(2,\zz)$-equivalence.
Since $GL(2,\zz)$ adds only the action of 
$ \left[{1 \atop 0} {0 \atop {-1}}\right]$,
it follows that all are equivalent to a form
with $B \ge 0$, and all of these are $GL(2,\zz)$-
inequivalent. 

Combining these two steps, where one takes
 $M=n^2$ in the second step, associates to 
any Descartes quadruple $(-n, x, y, z)$ the definite binary quadratic
form
\beql{N203}
[A, B, C] :=[-n+x, \frac{1}{2}(-n + x + y - z), -n + y]
\eeq
of discriminant $D :=(2B)^2 -4AC= -4n^2$. Furthermore
this form is reduced in the $GL(2,\zz)$- sense \eqn{N201} if and
only if $(-n, x,y,z)$ is a root quadruple. 
The primitivity condition \\ 
$gcd(-n, x,y,z) = 1$ transforms to 
$gcd(A,B, C) = 1$; conversely, $gcd(A,B, C) = 1$
implies $gcd(-n, x,y,z) = 1.$  ~~~$\Box$

Using Dirichlet's class number formula we obtain an explicit
formula for $N_{root}^{\ast}(-n)$.

\begin{theorem}~\label{Nth42}
For $n >1$ the number $ N_{root}^*(-n)$ of primitive integral root quadruples
satisfies
\beql{N204}
 N_{root}^*(-n)= \frac{n}{4} \prod_{p|n}(1-\frac{\chi_{-4}(p)}{p})+
2^{\omega(n)-\delta_n-1},
\eeq
where $\chi_{-4}(n) = (-1)^{(n-1)/2}$ for $n$ and $0$ for even $n$,
$\omega(n)$ denotes  the number of distinct primes dividing $n$, and 
$\delta_n =1$ if $n \equiv 2(\bmod~4)$ and $\delta_n=0$
otherwise.
\end{theorem}

\noindent{\bf Proof.}
The condition for a positive definite form $[A,B,C]$ to be 
$SL(2,\zz)$-reduced is that $0 \le 2|B| \le A \le C.$ 
The form classes under $GL(2,\zz)$-equivalence are derived from
the $SL(2, \zz)$-equivalence classes
using the coset decomposition
$$
GL(2, \zz) = \{ I, \left[{1 \atop 0} {0 \atop {-1}}\right] \} SL(2, \zz).
$$
The action of $$\left[{1 \atop 0} {0 \atop {-1}}\right]$$ maps
the set of $SL(2, \zz)$-reduced forms into itself, hence
\beql{N205}
h^{\pm}(-4n^2) = \frac{1}{2} (h(-4m^2) + a(-4m^2))
\eeq
where $a(-4m^2)$ denotes the number of primitive reduced ambiguous forms
of discriminant $-4m^2$, where we define a {\em reduced ambiguous form} to
be any $SL(2, \zz)$-reduced form of the shape
$[A, 0, C]$, $[B, 2B, C]$, and $[A, 2B, A]$,
with $B > 0.$ 
(These are exactly the reduced forms with $B \ge 0$
and $[A, B, C] \approx [A, -B, C]$,
see Mathews \cite[pp. 69--72 ]{Ma61}.)
A formula for the  general class number $h(d)$ was given
by Dirichlet, and we use the version in Landau~\cite{La27}.
For a discriminant $d<0$, Landau~\cite[Theorem 209]{La27},
states that
$$
h(d) = \frac{w_{d}}{2\pi} \sqrt{|d|}L(1, \chi_{d})
$$
in which  $L(s,  \chi_{d}) = \sum_{n=1}^\infty \chi_{d}(n) n^{-s}$,
with 
$\chi_{d}(n) = (\frac{d}{n})$ 
a real Dirichlet character $(\bmod~|d|)$,
and $w_{d}= 2$ if $d<-4$, with $w_{-4}=4$. Thus
for $d= -4n^2$ with $n > 1$, 
$$
h(-4n^2) = \frac{2n}{\pi} L(1, \chi_{-4n^2}).
$$ 
Landau~\cite[Theorems 214]{La27} also gives
$$
L(1, \chi_{-4n^2})= 
\prod_{p|n}( 1 -  \frac{\chi_{-4}(p)}{p}) ~L(1, \chi_{-4}).
$$
Combining this with
Lambert's formula $L(1, \chi_{-4}) = \frac{\pi}{4}$, we obtain,
$h(-4)=1$ and, for $n >1$, the class number formula
\beql{N206}
h(-4n^2) = \frac{n}{2} \prod_{p|n}( 1 -  \frac{ \chi_{-4}(p)}{p}).
\eeq

It remains to determine $a(-4n^2).$ We claim that, for $n > 1$, 
\beql{N206a}
a(-4n^2) = \left\{ \begin{array}{ll}
2^{\omega(n)-1} & \mbox{if}~~n \equiv 2~(\bmod~4), \\
2^{\omega(n)}   &  \mbox{otherwise.}
\end{array}
\right.
\eeq
To prove this we consider
for each $j \ge 0$ separately 
 the  cases $n = 2^j s$, with
$s$ odd. There are
three types of ambiguous reduced forms to consider:

(1) $[A, 0, C]$ with $0 < A \le C$; 

(2) $[2A, A, C]$ with $0 < 2A \le C$;

(3) $[A, B, A]$, with $0 < 2B \le A.$

\noindent The primitivity requirement is that $gcd(A,B,C)=1$.
For type (1), the discriminant condition gives
$ n^2=AC.$ Since $gcd(A,C)=1$ each relatively prime factorization
$n=pq$ with $gcd(p,q)=1$ gives a candidate pair $(A, C)=(p^2, q^2)$, 
and the requirement
$A\le C$ rules out half of these, provided $n > 1$, for
then $A=C$ cannot occur. We conclude that there are
exactly  $2^{\omega(n)-1}$ 
solutions of this type, when $n > 1.$
The  case $n=1$ is exceptional, with $ A= C = 1$ as the only
solution.

For the remaining count we must group types (2) and (3)
together.
For type (2) we have $-4n^2= 4A^2-8AC$ which gives
$n^2=A(2C-A)$. The condition $gcd(A, C)=1$ implies $gcd(A, 2C-A) = 1$ or $2$.
For type (3) we have $-4n^2= 4B^2 - 4A^2$ which gives
$n^2=(A+B)(A-B)$. The condition $gcd(A,B)=1$ 
implies $gcd(A+B, A-B) = 1$ or $2$.

We first suppose $j=1$. 
This is the easiest case because 
 there are no solutions of types (2) and (3).
Indeed, for type (2) we have $A(2C-A)$ is even, so $A$ must be even, and
$2C-A$ is even. But $gcd(A, C)=1$ means $C$ is odd, so $A$ and $2C-A$
are not congruent $(\bmod~4)$. Thus one of them is divisible by $4$,
hence $8$ divides $A(2C-A)=n^2$ which contradicts $n^2$ being divisible
by $4$ but not $8$. For type (3) we have $(A+B)(A-B)$ is even, so 
both factors are even. At least one of $A$, $B$ is odd since $gcd(A,B)=1$,
so the two factors are incongruent $(\bmod~4)$, hence one of them is
divisible by $4$, so again we deduce that $8$ divides $n^2$, a contradiction.

We next suppose $j=0$. Then  $n$ is odd and the primitivity requirement
for type (2) forms becomes  $gcd(A, 2C-A)=1$,
and for  type (3) forms becomes  $gcd(A+B, A-B)=1$. 
Consider a factorization $n^2=p^2 q^2$
with $gcd(p,q)=1$. For type (2), setting  this equal to the factorization
$A(2C-A)$ yields  $A=p^2, C=1/2(p^2 + q^2)$.
The condition $2A \le C$ becomes 
$$ p^2 \le \frac{1}{3}q^2.$$
For type (3),  setting $n=q^2p^2$ equal to the factorization
$(A+B)(A-B)$ yields $A= 1/2(q^2 + p^2),~B= 1/2(q^2 - p^2)$.
The condition $B \ge 0$ gives $p^2 \le q^2$ and the
condition $2B \le A$ gives $3p^2 \ge q^2$; thus
$ q^2/3 \le p^2 \le q^2$, and we have $ q^2/3 \ne p^2$ by the
relative primality condition.    
These two cases both
require $p^2 \le q^2$, and they have disjoint ranges for $p^2$ and
cover the whole range $0< p^2 \le q^2$. Thus exactly half the factorizations
lead to one solution and the other half to no solution,
and we obtain exactly $2^{\omega(n)-1}$ solutions of
types (2) and (3).

Finally we suppose  $j \ge 2$. In these cases 
the primitivity requirement for
type (2) forms becomes $gcd(A, 2C-A)=2$, and for 
 type (3) forms becomes  $gcd(A+B, A-B)=2$. We have
$n^2= 2^{2j}p^2q^2$ and there are several splittings of factors
to consider, for example $A=2^{2j-1}p^2, 2C-A=2q^2$ or
$A= 2q^2, 2C-A= 2^{2j-1}p^2$ in type (2) and 
$A+B= 2^{2j-1}q^2, A-B=2p^2$ or vice versa in case (3).
In all cases, a careful counting
of factorizations finds exactly half the possible
factorizations lead to solutions of  types (2) and (3), giving
$2^{\omega(n)-1}$ solutions of types (2) and (3). We omit
the details.

Adding up these cases gives $2^{\omega(n)}$ solutions in total when
$j=0$ or $j \ge 2$ and $2^{\omega(n)-1}$ solutions if $j=1$.
This proves the claim.

Combining \eqn{N205} and \eqn{N206} with the claim yields the
desired result.~~~$\Box$ \\

The exact formula for $N_{root}^*(-n)$  for prime $n=p$
gives $N_{root}^*(-p) = \frac{p+3}{4}$ if $p \equiv ~1(\bmod~4)$
and $N_{root}^*(-p) = \frac{p+5}{4}$ if $p \equiv ~3(\bmod~4)$,
as  in Table 2.
The exact formula also yields the following 
asymptotic upper and
lower bounds for $N_{root}^*(-n).$ 

\begin{theorem}~\label{Nth43}
There are positive constants $C_1$ and $C_2$ such that, for
all $n \ge 3$, 
\beql{N207}
C_1 \frac{n}{\log \log n} < N_{root}^{*}(-n) <  C_2~n \log \log n.
\eeq
\end{theorem} 

\noindent{\bf Proof.}
We use the exact formula \eqn{N204} of Theorem~\ref{Nth42}. The second term
$\frac{1}{2}a(-4n^2)$ grows like $O(n^\epsilon)$ and is
negligible compared to the first term. 
To estimate the first term we use inequalities of
the form
\beql{N208}
\prod_{p | n} (1 - \frac{1}{p}) \ge c_1 \frac{1}{\log \log n}.
\eeq
and
\beql{N209}
\prod_{p | n} (1 + \frac{1}{p}) \le c_2 \log \log n.
\eeq
To establish \eqn{N209}, we note that
its left side is upper bounded by  the product
taken over  the first $k$ distinct primes,
with $k$ chosen minimal so that  the product exceeds $n$.
Thus
$$
\log (\prod_{p | n} (1 + \frac{1}{p})) \le 
\sum_{p \le \log n / \log\log n} \log(1 + \frac{1}{p}) \le
\sum_{p \le \log n / \log\log n} \frac{1}{p} \le \log \log \log n 
$$
for large enough $n$, using \cite[Theorem 427]{HW}, and
exponentiating this yields \eqn{N209}.
The inequality \eqn{N208} follows from this because
$\prod_{p}(1 - \frac{1}{p^2})$ converges. 
~~~$\Box$

The bounds of Theorem~\ref{Nth43} are sharp up to a multiplicative constant.
We can take $n$ to be a product of all small primes $p \equiv 3~(\bmod~4)$
to achieve a sequence of values with
$N_{root}^{*}(-n) >> n \log\log n$, and to be a product of 
all small  primes $p \equiv 1~(\bmod~4)$  to achieve a sequence of values with
$N_{root}^{*}(-n) << \frac {n}{\log \log n}.$

\noindent\paragraph{Remark.} 
Theorem~\ref{Nth41} gives a one-to-one correspondence
between Descartes quadruples containing $-n$ and integral binary
quadratic forms of discriminant $-4n^2$ having nonnegative
middle coefficient. 
These Descartes configurations
are exactly the ones  that contain the  outer circle  of the packing
as one of their circles.
Under this correspondence
the reduction algorithm of \S3 for Descartes quadruples 
matches the Gaussian reduction  algorithm for positive definite binary
quadratic forms. 

A related  correspondence holds 
for the exceptional integral
Apollonian packing $(0,0,1,1)$ having $n=0$.
There is a unique 
primitive reduced positive definite form of determinant $0$, namely
$Q= [0, 0, 1] = U^2$, which is an 
ambiguous form.
In this case the special Descartes quadruples 
%having a binary form interpretation 
are those that contain  one of the two infinite
edges of the packing as one of their circles. 
They are essentially Descartes configurations of ``Ford circles''.
and on these configurations the reduction algorithm corresponds
to an additive variant of the continued fraction
algorithm. For futher information see
 Ford~\cite{Fo38} or  Rademacher~\cite[pp. 41--46]{Ra64}, 
\cite[pp. 264--267]{Ra73}.

%--------------------------------------------------------------------
%
% Section 5
%
%--------------------------------------------------------------------
%

\section{Integers Represented by a Packing: Asymptotics}
\setcounter{equation}{0}
In this section we  
study the ensemble of integer curvatures that occur in an
integer Apollonian circle packing, where integers are counted with
the multiplicity that they occur in the packing. Their asymptotics
are known to be related to the Hausdorff dimension of the residual
set of the packing, as follows from work of Boyd described in
Theorem~\ref{th72} below. At the end of the section we 
begin the study of the set of integer curvatures that occur, counted
without multiplicity.

The {\em residual set} of a disk packing $\sP$ 
(not necessarily an Apollonian packing) is the set
remaining after all the (open) disks in the packing are
removed, including any disks with ``center at infinity.''
For a general disk packing $\sP$,
we denote the 
Hausdorff dimension of the residual set by $\alpha(\sP)$ and
call it the {\em residual set dimension} of the packing.
The definition of Hausdorff dimension
can be found in Falconer \cite{Fa86}, who also studies the residual
sets of Apollonian packings in \cite[pp. 125--131.]{Fa86}. 

The residual sets of Apollonian packings all have
the same Hausdorff dimension, which we denote by
$\alpha.$ This is a consequence of the equivalence
of such residual sets under M\"{o}bius  transformations
(see \cite[Sect. 2]{GLMWY21}),
using also the fact that the Hausdorff dimension 
strictly exceeds one, as follows from results described below.

The {\em exponent} or {\em packing constant} $e(\sP)$ of a 
bounded circle packing $\sP$ 
(not necessarily an Apollonian packing) is defined to be
\[
e(\sP) :=\sup\{ e: \sum_{C \in \sP}  r(C)^e=\infty\}=
\inf\{e : \sum_{C\in \sP}  r(C)^e < \infty\},
\]
in which  $r(C)$ denotes the radius of the circle $C$.
This number has been extensively studied in the literature, beginning
in 1966 with the work of Melzak \cite[Theorem 3]{Me66}, who showed
that in any circle packing that covers all but  
a set of measure zero one has  
$\sum_ {C \in \sP}  r(C) =\infty$. 
He constructed a circle packing with $e(\sP) = 2$ and 
showed for Apollonian packings that
$e(\sP)$ lies strictly between $1.035$ and $1.99971$. 
He conjectured that the minimal value of
$e(\sP)$ is attained by an Apollonian  circle packing.
In 1967 J. Wilker~\cite{Wi67} showed that all osculatory circle
packings $\sP$, which include all Apollonian circle packings, have
the same exponent $e(\sP)$, which we call the 
{\em osculatory packing exponent e.}
He also showed that $e \geq 1.059.$
Later Boyd \cite{Bo71}, \cite{Bo73a}, \cite{Bo82}
improved this to  $1.300 < e < 1.314$.
 Recent non-rigorous computations 
of Thomas and Dhar~\cite{TD94} estimate the Apollonian packing exponent
to be $1.30568673$ with a possible error of $1$ in the last digit.

The relation between the packing exponent and the 
residual set dimension of Apollonian packings was resolved
by an elegant result of D. Boyd~ \cite{Boyd73}.

\begin{theorem}~\label{th71}
{\rm (Boyd)} The exponent $e$ of any bounded
Apollonian circle packing is equal to
the Hausdorff dimension $\alpha$ of the residual set of any
Apollonian circle packing.
\end{theorem}

The inequality $e \geq \alpha$ follows from a 1966 result of
Larman~\cite{Lar66}, and in 1973 Boyd proved the 
matching upper bound $ \alpha \geq e.$
A simpler proof of the upper bound was later given by C. Tricot~\cite{Tri84}.

Given a bounded circle packing $\sP$ we
 define the {\em circle-counting function}  $N_{\sP}(T)$ 
to count the number of circles in the packing
whose radius of curvature is no larger than $T$, i.e., whose radius is at
least $\frac {1}{T}.$ Boyd~\cite{Bo82} proved the following improvement
of the result above.

\begin{theorem}~\label{th72}
{\rm (Boyd)} For a bounded Apollonian circle packing $\sP$, 
the circle-counting
function $N_{\sP}(T)$ satisfies
\beql{N500}
 \lim_{T \to \infty} \frac {log~ N_{\sP}(T)} {log~ T} = \alpha,
\eeq
where $\alpha$ is the Hausdorff dimension of the residual set. That is,
$N_{\sP}(T) = T^{\alpha + o(1)}$ as $T \to \infty.$
\end{theorem}.
Theorem \ref{th42b} showed that the curvatures of all circles
in the packing, excluding the root quadruple,
can be enumerated by the elements of the Apollonian group $\sA$.
>From this one can derive a relation between
 the number of elements of $\sA$ having height below 
a given bound $T$ and the Hausdorff dimension $\alpha$.
We measure the {\em height} of an element $M \in \sA$ using the Frobenius norm
\beql{N501}
||M||_F := (tr[ M^T M])^{1/2} = (\sum_{i, j} M_{ij}^2)^{1/2}.
\eeq

\begin{theorem}~\label{th73}
The number of elements $N_T(\sA)$ of height at most $T$ in the 
Apollonian group
$\sA$ satisfies
\beql{N502}
N_T(\sA) = T^{\alpha + o(1)},
\eeq
as $T \to \infty,$
where $\alpha$ is the Hausdorff
dimension of the residual set of any Apollonian packing.
\end{theorem}

In order to prove this result, we establish two preliminary lemmas.

\begin{lemma}\label{le51a}
Let $M= S_{i_m} \cdots S_{i_2} S_{i_1} \in \sA$, the Apollonian group,
and suppose that $i_j \neq i_{j+1}$ for $1 \le j \le m-1$, and $m \ge 2$.
In each row $k$ of $M$,
\begin{itemize}
\item[(i)]
$M_{kl} \le 0$ if $l =i_1$,
\item[(ii)]
$M_{kj} \ge |M_{kl} |$ for $l=i_1$ and $j \neq l$.
\end{itemize}
\end{lemma}

\pf The lemma follows by induction on $m$.
It is true for $m=1$, since each matrix $S_i$ has $i^{\rm th}$ 
column negative (or zero).

Suppose (i)--(ii) hold for $M' = S_{i_m} \cdots S_{i_2}$.
Suppose, for convenience, that $i_1 =1$. Then
$$
M = M' S_{i_1} =
\left[
\begin{array}{cccc}
-M'_{11} & 2M'_{11} + M'_{12} & 2M'_{11} + M'_{13} & 2M'_{11} + M'_{14} \\
-M'_{21} & 2M'_{21} + M'_{22} & 2M'_{21} + M'_{23} & 2M'_{21} + M'_{24} \\
-M'_{31} & 2M'_{31} + M'_{32} & 2M'_{31} + M'_{33} & 2M'_{31} + M'_{34} \\
-M'_{41} & 2M'_{41} + M'_{42} & 2M'_{41} + M'_{43} & 2M'_{41} + M'_{44} \\
\end{array}
\right] \,.
$$
Since $i_2 \neq i_1 =1$ all $M'_{i1} \ge 0$ by (ii) 
of the induction hypothesis, so
$M_{i1} = M'_{i1} \le 0$ gives (i).
Next, note that
$$M_{kj} = 2M'_{k1} + M'_{kj}
\ge 2M'_{k1} - | M'_{kl} | \ge M'_{k1} = | M_{k1} |
$$
since $M'_{kj} \ge |M'_{kj} |$ and $M'_{kj} \ge - | M'_{kl} |$ in all cases
by (ii).
This completes the induction step in this case. The arguments
when $i_1 = 2, 3,$ or $4$ are similar.~~~$\Box$

\begin{lemma}\label{le52a}
Let $\bv = (a,b,c,d)^T$ be an integer root quadruple with $a< 0$.
Then there are positive constants $c_0 = c_0 (\bv )$ and 
$c_1 = c_1 ( \bv )$ depending on $\bv$ such that
\beql{DI52}
c_0 \| M \|_F \le ||M \bv ||_\infty \le c_1 \| M\|_F , \quad\mbox{for~all}\quad
M \in \sA ~.
\eeq
\end{lemma}

\pf
For the upper bound, we have
\beql{DI53}
|| M\bv ||_\infty
\le 2 |M \bv | \le 2 \| M \|_F | \bv | ~,
\eeq
so we may take $c_1 = 2 | \bv |$.

For the lower bound, we first show that if
$M = S_{i_m} \cdots S_{i_2} S_{i_1}$ with $i_j \neq i_{j+1}$ and $i_1 =1$,
we have
\beql{DI54}   
|| M \bv ||_\infty \ge \frac{1}{2} \| M \|_F ~.
\eeq
The vector $\bv$ has sign pattern $(-, +,+,+ )$ and Lemma \ref{le51a} 
shows that $M$ has first column
nonpositive elements and other columns nonnegative.
Thus all terms in the product $M \bv$ are nonnegative, and hence
$$(M \bv )_i \ge \sum_{j=1}^4 | M_{ij} | | \bv_j | \ge
\sum_{j=1}^4 |M_{ij} | ~,
$$
because $a< 0$ implies
$\min ( |a|, |b| , |c|, |d| ) \ge 1$.
Thus
$$||M \bv ||_\infty \ge \frac{1}{4}
\sum_{i,j} | M_{ij} | \ge \frac{1}{2} \| M\|_F ~.
$$

It remains to deal with the cases where $i_1 = 2,~3$ or $4.$ 
By Theorem \ref{th42b}, the value $||M \bv||_\infty$
gives the curvature of a particular circle in the packing, and this circle 
lies in one of the four lunes pictured in Figure 3 according to 
the value of $i_m$.

\begin{figure}[htbp]\label{fig3}
\centerline{\epsfxsize=2.5in \epsfbox{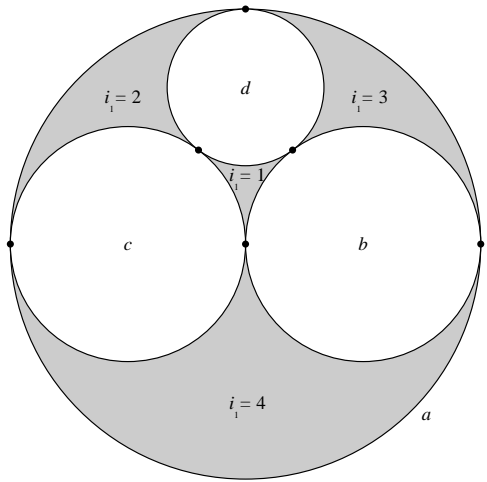}}
\caption{Four lunes of Descartes quadruple.}
\end{figure}

The bound \eqn{DI54} applies to all circles in the central lune 
corresponding to $i_m =1$.
For the remaining
cases, we use the fact that there exists a M\"{o}bius transformation 
$\phi: \hat{\cc} \to \hat{\cc}$ with $\phi \in Aut ( \sP )$, which 
fixes the Descartes configuration corresponding to $\bv$ but cyclically 
permutes the four circles $a \to b \to c \to d$.
In particular $\phi$ also cyclically permutes the four 
lunes $i_1 =1 \to i_1 =2 \to i_1 =3 \to i_1 =4$.
Now $\phi$ maps the center of circle $d$ to the center of circle 
$a$, which is the point at
infinity, and maps the point at infinity to the center of circle $b$.
It follows that the stretching
factor of the map $\phi$ inside the four lunes is bounded above and below 
by positive absolute constants $c_2$ and $c_2^{-1}$.
Since $\phi$ maps the lune $i_1=4$ to $i_1=1$ we conclude for cases where
$i_1=4$ that
\beql{DI55}
|| M \bv ||_\infty \ge \frac{1}{2c_2} \| M \|_F ~.
\eeq
Applying the same argument
to $\phi^2$ and $\phi^3$ gives the similar bound for the cases 
$i_1=3$ and $i_1=2$.
We conclude that the lower bound in \eqn{DI52} holds with 
$c_0 = \frac{1}{2c_2}$.~~~$\Box$

\paragraph{Proof of Theorem \ref{th73}.}
Pick a fixed quadruple having $a< 0$, say $\bv = (-1,2,2,3)$, and let
$\sP_{\bv }$ be the associated Apollonian packing.
By Theorem \ref{th42b}, each $M \in \sA$ corresponds to a circle of 
curvature $||M \bv ||_\infty$ in $\sP_{\bv }$, and all circles are so labelled 
except the four circles in $\bv$.
Lemma \ref{le52a} shows that each $\| M\|_F < T$ produces a circle of
curvature at most $c_1 T$.
Now Theorem \ref{th72} asserts there are at most $T^{\alpha + o(1)}$ 
such circles, hence $N_T ( \sA ) \le T^{\alpha + o(1)}$.
Conversely Lemma \ref{le52a} implies that each circle of curvature 
$|| M \bv||_\infty \le T$ comes from a matrix $M \in \sA$ with 
$\| M \|_F \le \frac{1}{c_0} T$.
Since there are at least $T^{\alpha + o(1)}$ such circles, we obtain
$N_T ( \sA ) \ge T^{\alpha +o(1)}$, as desired.~~~$\Box$ 

Can the estimate of Theorem~\ref{th73} be sharpened to obtain an asymptotic
formula? A. Gamburd has pointed out to us that the method of 
Lax and Phillips~\cite{LP82} might prove useful in studying this question.

We now turn to a different question: How many
different integers occur, counted without multiplicity, in a given
integral Apollonian circle packing $\sP_{\bv }$?
This seems to be a difficult problem.
It is easy to prove that at least $cT^{1/2}$ of all
integers less than $T$ occur in a given packing. This comes
from considering the largest elements of
the vectors $\{(S_1 S_2)^j \bv: j = 1, 2,...\},$
where $\bv$ is a root quadruple, which are curvatures in the
packing, by Theorem ~\ref{th42b} above. These values grow like
$j^2$ (see the example (1) in \S7). Concerning  the true
answer to the question above, we propose the
following conjecture.

\noindent {\bf Positive Density Conjecture}.
{\em Each integral Apollonian packing represents a positive fraction
of all integers. }

Theorem~\ref{th72}
shows that 
the average number of representations
of an integer $n$ grows like  $n^{\alpha-1}$, which goes rapidly to 
infinity as $n \rightarrow \infty$. Therefore one might guess
that all sufficiently large integers are represented. However
in the next section we will show there are always some  congruence
restrictions on which integers occur. There we
formulate  a stronger version of this conjecture and present
numerical evidence concerning it. 

%--------------------------------------------------------------------
%
% Section 6
%
%--------------------------------------------------------------------
%

\section{Integers Represented by a  Packing: Congruence Conditions}
\setcounter{equation}{0}
In this section we study congruence restrictions on
the set of integer curvatures which occur in a 
primitive integral Apollonian packing. 

We first show that there are always congruence restrictions
(mod~ 12).

\begin{theorem}~\label{th61}
In any primitive integral Apollonian packing, the 
Descartes quadruples
(mod~12) all fall in exactly one of four possible orbits. 
The first orbit $Y$ (mod~12) consists of all permutations
of 
\beql{601a}
\{ (0,0,1,1),~(0,1,1,4), (0,1,4,9), (1,4,4,9), (4,4,9,9)\}~(\bmod~12).
\eeq
The other three orbits are 
$(3,3,3,3) - Y$, $(6,6,6,6) + Y$ and $(9,9,9,9) - Y$~ (mod~12).
\end{theorem}

\paragraph{Remark.} Each orbit contains only $4$ different residue
classes (mod~12), hence $8$ residue classes (mod~12) are excluded 
as curvature values.

{\bf Proof.} A straightforward computation, using the action of
the Apollonian group (mod ~12), shows that (up to ordering) the set of all
quadruples without common factors of $2$ or $3$ (mod $12$)
consists of the list below, which are grouped into eight orbits
under the action of the Apollonian group (mod~12). 

\begin{eqnarray*}
\begin{array}{ccclllll}
 (1)& Y &=&(0,0,1,1)   &(0,1,1,4)  &(0,1,4,9) &(1,4,4,9) &(4,4,9,9); \\
 (2)&(3,3,3,3)-Y&=&(6,6,11,11) &(2,6,11,11)&(2,3,6,11)&(2,2,3,11)&(2,2,3,3); \\
 (3)&(6,6,6,6)+Y&=&(3,3,10,10) &(3,6,7,10)&(3,7,10,10)&(6,7,7,10)&(6,6,7,7); \\
 (4)&(9,9,9,9)-Y&=&(0,0,5,5)   &(0,5,5,8)  &(0,5,8,9) &(5,8,8,9) &(8,8,9,9); \\ 
 (5)& -Y&=&(0,0,11,11) &(0,8,11,11)&(0,3,8,11)&(3,8,8,11)&(3,3,8,8); \\
 (6)&(3,3,3,3)+Y&=&(0,0,7,7)   &(0,4,7,7)  &(0,3,4,7) &(3,4,4,7) &(3,3,4,4); \\
 (7)&(6,6,6,6)-Y&=&(2,2,9,9)   &(2,2,5,9)  &(2,5,6,9) &(2,5,5,6) &(5,5,6,6); \\
 (8)&(9,9,9,9)+Y&=&(9,9,10,10) &(1,9,10,10)&(1,6,9,10)&(1,1,6,10)&(1,1,6,6); \\
\end{array} 
%\ ({\rm mod }~12)
\end{eqnarray*}

To check the orbit structure is as given, note that the action
of the four generators
of the Apollonian group on the five elements of 
the orbit $Y$ is summarized in
the following transition matrix:
\[
\frac{1}{4} \left( \begin{array}{ccccc}
                   2 & 2 & 0 & 0 & 0 \\
                   1 & 1 & 2 & 0 & 0 \\   
                   0 & 1 & 2 & 1 & 0 \\
                   0 & 0 & 2 & 1 & 1 \\
                   0 & 0 & 0 & 2 & 2 \end{array}\right).
\]
We may  view this matrix as the
transition matrix of 
a Markov chain (after rescaling each row to be stochastic),
and find that the action is transitive and 
the stationary distribution is
$(\frac{1}{10}, \frac{1}{5}, 
\frac{2}{5}, \frac{1}{5},  \frac{1}{10})$.
The other seven orbits
have the same transition 
matrix and the same stationary distribution as $Y$. 

There exist integral solutions to the
Descartes equation in all of the congruence classes (mod~12) in
the list above. However we
recall that a Descartes quadruple $(a, b, c, d)$ coming from
an Apollonian packing satisfies the extra condition
\beql{N601}
a + b + c + d > 0.
\eeq
In the rest of the proof
we show that this extra condition excludes half of the orbits above,
namely orbits (5)- (8).

As a preliminary, we observe that any integer solution
$(a, b, c, d)$ to the Descartes equation \eqn{descartes} yields
a unique integer solution to the equation
\beql{N602}
4m^2 + 4a^2 + n^2 = l^2,
\eeq
and vice-versa. Here the solution to \eqn{N602} is
given by 
\beql{N603}
\left[\begin{array}{c}
a \\ n \\ l \\ m
\end{array} \right] =
\left[
\begin{array}{rccr}
1 & 0 & 0 & 0 \\
0 & 1 & -1 & 0 \\
2 & 1 & 1 & 0 \\
 -\frac{1}{2}  &  -\frac{1}{2}  & -\frac{1}{2} & \frac{1}{2}
\end{array}
\right]
\left[\begin{array}{c}
a \\ b \\ c \\ d
\end{array} \right] =
\left[ \begin{array}{c}
a \\ b - c \\ 2a + b + c \\ \frac{1}{2}(d - a - b  - c)
\end{array}
\right] \,.
\eeq
In the reverse direction, an integer solution to \eqn{N602} gives
one to the Descartes equation via 
\beql{N604}
\left[\begin{array}{c}
a \\ b \\ c \\ d
\end{array} \right] =
\left[
\begin{array}{rccr}
1 & 0 & 0 & 0 \\
-1 & \frac{1}{2} &  \frac{1}{2}    & 0 \\
-1 & -\frac{1}{2} &  \frac{1}{2} & 0 \\
-1  &  0  & 1 & 2
\end{array}
\right]
\left[\begin{array}{c}
a \\ n \\ l \\ m
\end{array} \right] =
\left[ \begin{array}{c}
a \\ \frac{1}{2}( l - 2a + n) \\ \frac{1}{2}( l - 2a - n) \\ 
2m + l - a
\end{array}
\right] \,.
\eeq
Solutions to the Descartes equation satisfy a congruence $(\bmod ~2)$
which guarantee that the maps above take integral solutions to 
integral solutions, in both directions.
Now \eqn{N602} gives
\beql{N604a}
 l^2 \geq 4a^2 + m^2 \geq 2(|a| + |m|)^2 \geq (|a| + |m|)^2, 
\eeq
and equality holds if and only if $\ell = a = m = 0.$ 
%Using   \eqn{N602} this implies $n = 0$, hence
%and by \eqn{N604}
%this gives $a = b = c = d = 0.$ 
In particular, if $\ell > 0$, then \eqn{N604a} gives
\beql{N604b}
\ell > |a| + |m|.
\eeq 

We assert that any  integer
 solution $(a, b, c, d)$  to the Descartes equation has
\beql{N605}
a + b + c + d > 0 \qquad\mbox{if~ and~ only~ if}~~~~ l > 0.
\eeq
To prove this, note that
if $a + b + c + d > 0$ then by Lemma~\ref{le31a} (i) we have 
$$\ell = 2a + b + c = (a + b) + (a + c) \geq 0.$$
Equality can  hold here only if $a= b= c= 0,$ which implies
$d = 0$, which contradicts the assumption $a + b + c + d > 0$ .
Conversely, if $\ell > 0$, then, using \eqn{N604b},
$$a + b + c + d = 2\ell + 2m - 2a \geq 2(\ell - |a| - |m|) > 0,$$
so \eqn{N605} is proved.

{\bf Claim :} No primitive integer
Descartes quadruples with $a +  b + c + d > 0$ occur in the orbits (5)--(8). 
 
We prove the claim for orbit (8);  the arguments to rule out
orbits (5), (6), (7) are similar.
We argue by contradiction. Suppose there were such a solution
in orbit (8). Since the Apollonian group acts transitively
on the orbit, and preserves the condition $a + b + c + d > 0,$
there would be such a quadruple 
$(a,b,c, d) \equiv (1, 1, 6, 6) (\bmod~ 12).$
In this case $l = 2a + b + c \equiv 9 (\bmod~ 12)$ and
$$m \equiv \frac{1}{2}(6 - 6 + 1 + 1) \equiv 1 (\bmod~ 6),$$
which gives $m^2 \equiv 1(\bmod~ 12).$ Now \eqn{N602} gives
\beql{N606}
 (l + 2m) (l - 2m) = l^2 - 4m^2 = 4a^2 + n^2 > 0. 
\eeq
Since $a + b + c + d >0$ we
 have $l > 0$ by \eqn{N605}. Then in the equation above at 
least one
of the factors on the left side must be  positive, hence they both are.
Consider $l + 2m > 0.$ We have
\beql{N607}
l + 2m \equiv 9 \pm 2 (\bmod~12) \equiv 3 (\bmod 4).
\eeq
Consider any prime $p \equiv 3 (\bmod ~ 4)$ dividing $l + 2m$.
Then it divides $4a^2 + n^2$, which it must divide to an
even power, say $p^{2e}$, with $a \equiv n \equiv 0 (\bmod~p^e).$
If $p$ also divides $l - 2m$, then it would divide both $l$ and $m$,
and then \eqn{N604} would imply that it divides $\gcd(a,~b,~c,~d)$,
which contradicts the primitivity assumption $\gcd(a,~b,~c,~d)= 1$.
Therefore $p$ does not divide $ l - 2m,$ and we conclude from \eqn{N606}
that $p^{2e} || l + 2m$. It follows that all primes $p \equiv 3 (\bmod ~4)$
that divide $l + 2m$ do  so to an even power, hence we must have
have $l+ 2m \equiv 1 (\bmod~4),$ a contradiction. This rules out
orbit (8), which proves the claim in this case. 

Theorem~\ref{th61} follows from the claim.
 $\Box$

At the end of this section we
present numerical evidence that suggests 
that these congruences (mod $12$)
are the only congruence restrictions
for the integer  packing   $(-1, 2, 2, 3)$.
However
there are stronger modular restrictions (mod $24$) that
apply to  other integer packings.
For example, in the packing $(0,0,1,1)$ (Fig. 2), any curvature which
occurs must be congruent to $0,1,4,9,12$ or $16$ (mod $24$) (these are
the quadratic residues modulo $24$). Thus only $6$ classes (mod $24$)
can occur rather than the $8$ classes allowed by Theorem~\ref{th61}.

It seems likely that the full set of congruence restrictions possible
(mod $m$) is attained for $m$ a small fixed power $2^a 3^b$, perhaps
even $m=24.$ We are a long way from proving this. As evidence in
its favor, we prove the following result, which 
shows that all residue classes modulo $m$
do occur for any $m$ relatively prime to $30$.

\begin{theorem} Let $\mathcal P$ be a primitive integral 
Apollonian circle packing. For any integer $m$ with 
gcd$(m,30)=1$,  every residue class
modulo $m$ occurs as the value of some circle curvature in the
packing $\mathcal P$.
\end{theorem}
\pf
Observe that the $s$-term product $W(s)=\ldots S_2S_1S_2S_1$ is
$$\left( \begin{array} {cccc}
-s & s+1 & s(s+1) & s(s+1) \\
-(s-1) & s & s(s-1) & s(s-1) \\
0 & 0 & 1 & 0 \\
0 & 0 & 0 & 1 \end{array} \right)$$
(where the top two rows are interchanged if $s$ is even). Of course,
the two non-trivial rows can be placed anywhere by choosing the
two matrices from the set ${S_1,S_2,S_3,S_4}$ appropriately.

If $(a,b,c,d)^T$ is a quadruple in $\mathcal P$ then the product
\beql{eq6.1}
W(s)(a,b,c,d)^T=~~~~~~~~~~~~~~~~~~~~~~~~~~~~~~~~~~~~~~~~~~~~\nonumber\\
(-sa+(s+1)b+s(s+1)c+s(s+1)d,-(s-1)a+sb+s(s-1)c+s(s-1)d,c,d)^T
\eeq
is also in $\mathcal P$
as well. Let $\mathcal J$ denote the set of all rows 
$(\alpha,\beta,\gamma,\delta)$ which can occur in a product of matrices taken from
${S_1,S_2,S_3,S_4}$. Thus, if $(\alpha,\beta,\gamma,\delta) \in \mathcal J$
then so are:\\
$(-\alpha, 2\alpha + \beta, 2\alpha + \gamma, 2\alpha + \delta)$,
$(2\beta + \alpha, -\beta, 2\beta +\gamma, 2\beta + \delta)$,
$(2\gamma + \alpha, 2\gamma + \beta, -\gamma, 2\gamma + \delta)$,\\ and
$(2\delta + \alpha, 2\delta + \beta, 2\delta + \gamma, -\delta)$.
Therefore,\\
$(-\alpha, 2\alpha + \beta, 2\alpha + \gamma, 2\alpha + \delta),
(3\alpha+2\beta, -2\alpha-\beta,6\alpha+2\beta+\gamma, 6\alpha + 2\beta +\delta),\\
(-3\alpha -2\beta, 4\alpha + 3\beta, 12\alpha + 6\beta + \gamma, 12\alpha + 6\beta + \delta),\ldots$, and in general,
\beql{eq6.2}
(-r\alpha - (r-1)\beta, (r+1)\alpha + r\beta, r(r+1)\alpha + r(r-1)\beta + \gamma, r(r+1)\alpha + r(r-1)\beta + \delta)
\eeq
are all in $\mathcal J$ for all $r$
(as well as all permutations of these). Now
substitute $(\alpha, \beta, \gamma, \delta) = (s(s+1), s(s+1), -s, s+1)
\in \mathcal J$ into \eqref{eq6.2}. This shows that the row vector
\beql{eq6.3}
\rho := (-(2r-1)s(s+1), (2r+1)s(s+1), 2r^2s(s+1)-s, 2r^2s(s+1)+s+1) \in
\mathcal J.
\eeq
The sum of the last two coordinates of $\rho$ is
\beql{eq6.4}
4r^2s(s+1)+1 = r^2((2s+1)^2 - 1) +1 = r^2(x^2 - 1) + 1 = u^2-r^2+1
\eeq
where $u=rx$ and $x=2s+1$. The g.c.d. of these two summands
must divide their difference, which is $2s+1$. It is well known (and
easy to show) that for any prime power $p^w$ with $p > 5$, in at least
one of the pairs $\{1,2\}, \{4,5\}$ and $\{9,10\}$ are both 
nonzero quadratic
residues modulo $p^w$. For each $p~|~m$, let $\{a_p,a_p+1\}$ denote such
a pair. Define $u_p$ and $r_p$ so that
\beql{eq6.5}
u_p^2 \equiv a_p, \qquad r_p^2 \equiv a_p+1~(\mbox{mod}~ p^w)
\eeq
where $p^w$ is the largest power of $p$ dividing $m$. Since gcd$(r_p,p)=1$
then we can define $x_p \equiv u_p r_p^{-1}~(\mbox{mod}~p^w)$. We can
guarantee that $x_p$ is odd by adding a multiple of $p_w$ if necessary.
Hence, for these choices, the expression in \eqref{eq6.4} 
is $0$ modulo $p^w$, i.e.,
\beql{eq6.6}
r_p^2(x_p^2-1) + 1 \equiv 0~(\mbox{mod}~ p^w).
\eeq
Of course, we can use the values $r_p + kp^w$ and $x_p + lp^w$ in place
of $r_p$ and $x_p$ in \eqref{eq6.6} for any $k$ and $l$. Note that
gcd$(x_p,p) = 1$. Letting $p$ range over all prime divisors of $m$, then
by the Chinese Remainder Theorem, there exist $X$ (odd) and $R$ such that
\beql{eq6.7}
R^2(X^2-1) + 1 \equiv 0~(\mbox{mod}~ p^w)
\eeq
for all $p^w |~m$. Thus,
\beql{eq6.18}
R^2(X^2-1) + 1 \equiv 0~(\mbox{mod}~ m), ~~~~gcd(X,m) = 1.
\eeq
Hence, by \eqref{eq6.3} and \eqref{eq6.4} we can find a row 
modulo $m$ in $\mathcal J$
of the form $(C,D,A,-A) (\mbox{mod}~ m)$ where it easy to check 
that gcd$(A,m) = 1$.

We can now apply the transformation preceding \eqref{eq6.2}
to $(C,D,A,-A)$ to get
the following rows modulo $m$ in $\mathcal J$:
\begin{eqnarray*}
(C,&D,&A, -A)~~(\mbox{mod}~m)\\
(2A+C,&2A+D,& -A, A)~~(\mbox{mod}~m)\\
(4A+C,&4A+D,&A, -A)~~(\mbox{mod}~m)\\
(6A+C,&6A+D,&-A, A)~~(\mbox{mod}~m)\\
&\cdots&
\end{eqnarray*}
and more generally
\beql{eq6.8}
(4tA+C, 4tA+D, A, -A)  ~~(\mbox{mod}~m) \in \mathcal J ~~\mbox{for all}~t \geq 0.
\eeq
Suppose for the moment (and we will prove this shortly) that we can
find $(a,b,c,d)^T \in \mathcal P$ with gcd$(a+b,m) = 1$. Taking the
inner product of the row in \eqref{eq6.8} with $(a,b,c,d)^T$ , we get
the curvature value
\beql{eq6.9}
(4tA+C,4tA+D,A,-A) \cdot (a,b,c,d)^T&(\mbox{mod}~m)\\
\equiv 4A(a+b)t + Ca +Db + Ac -Ad&(\mbox{mod}~m).
\eeq
Since gcd$(4A(a+b),m) = 1$ then these values range over a complete
residue system modulo $m$ as $t$ runs over all positive integers.
 The proof will be complete now if
we can establish the following result.\\
{\bf Claim}.~
If $\mathcal P$ is a primitive packing then for any odd $m \geq 1$, there
exists $(a,b,c,d)^T \in \mathcal P$ with gcd$(a+b,m) = 1$.\\
{\bf Proof of Claim}.~
First recall that for any $(a,b,c,d) \in \mathcal P$, we have
gcd$(a,b,c) = 1$. We have also seen by \eqref{eq6.1}, if $(a,b,c,d) \in
\mathcal P$ then for any $r > 1$,
\beql{6.10}
(A(r),B(r),C(r),D(r)):=~~~~~~~~~~~~~~~~~~~~~~~~~~~~~~~~~~~~~~~~\nonumber\\
(-ra+(r+1)b+r(r+1)(c+d),-(r-1)a+rb+r(r-1)(c+d),c,d) \in \mathcal P
\eeq
as well. Define $q(r)$ to be the sum of the first two components of
this vector:
\[
q(r):=A(r)+B(r)=2(c+d)r^2 - 2(a-b)r +a+b.
\]
Let $p$ denote a fixed odd prime. We show that 
\beql{614a}
q(r) \not \equiv 0~(\mbox{mod}~p)~~\mbox{for~some}~~r \geq 1.
\eeq
Suppose to the contrary that 
$q(r) \equiv 0~(\mbox{mod}~p)$ for all $r$. Thus,
\begin{eqnarray*}
q(0) &\equiv& a+b \equiv 0~(\mbox{mod}~p)\\
q(1) &\equiv& 2(c+d)-2(a-b)+(a+b)
\equiv  2(c+d)-a+3b \equiv 0~(\mbox{mod}~p)\\
q(2) &\equiv& 8(c+d) - 4(a-b) + (a+b) \equiv 8(c+d)-3a+5b \equiv 0~
(\mbox{mod}~p)
\end{eqnarray*}
which implies $a \equiv b \equiv c+d \equiv 0~(\mbox{mod}~p).$
However, since\\
$a = b+c+d \pm 2 \sqrt{b(c+d)+c~d}$
then
$cd \equiv 0~(\mbox{mod}~p)$, i.e.,
$c \equiv 0$ or $d \equiv 0~(\mbox{mod}~p)$.\\
This would imply that  $\mathcal P$ is not primitive, a contradiction
which establishes \eqn{614a}.

To finish proving the claim, for each $p|m$
let $r_p$ satisfy $q(r_p) \not \equiv 0~(\mbox{mod}~p)$.
Then we have\\
$$q(r_p+kp) \equiv q(r_p) \not \equiv 0~(\mbox{mod}~p)$$
for all $k \geq 0.$ By the Chinese Remainder Theorem 
one can find  $R$ and $S$ such that 
$q(R+kS) \not \equiv 0~(\mbox{mod}~p)$ for all $p~|~m$ and all $k$. 
In particular,
$\gcd(q(R),m) = \gcd(A(R) + B(R), m) = 1$,
and the Claim is proved.
$\Box$

Which integers occur as curvatures, when the congruence conditions are
taken into account?
We consider numerical data for two cases. The first case is the packing
with root
quadruple  $(-1, 2, 2, 3)$, where Theorem~\ref{th61} permits only
values $2, 3, 6$, or $11$ (mod $12$).
Not all such integers   
appear in the Apollonian packing $(-1, 2, 2, 3)$,
for example in the class $6$ (mod $12$) the value 78 is missed. 
In Table 3 we present the missing values in these residue
classes for the first million
integers. Only 61
integers congruent to $2, 3$, or $6$ do not occur in the
packing $(-1, 2, 2,3)$, the largest being $97287$ (see Table 3),
and no integers 11 (mod $12$) are missed.
This data suggests that there
are finitely many missing values in total, 
with $97287$ being the largest one.

\begin{table}\label{mis}
\centerline{$n \equiv 3$ (mod 12)}
\begin{center}
\begin{tabular}{|llllllllll|}
\hline
159 & 207 & 243 & 435 & 603  & 711 & 1923  & 2175& 2319 &   3711 \\
4167& 4959& 4995& 5283& 6015 & 6879& 7863  & 10095&  10923& 11295 \\
12063 & 16311& 16515& 18051& 19815  &21135 & 23175& 28323& 41655&   48075 \\
68055 & 97287& & & & & & & & \\
\hline 
\end{tabular}
\end{center}
\vspace{.2cm}

\centerline{$n \equiv 6$ (mod 12)}
\begin{center}
\vspace{-.3cm}
\begin{tabular}{|llllllllll|}
\hline
78 & 246 & 342 & 834& 1422& 2010& 2022 & 2454 & 2718 & 2766 \\
3150 & 3402 & 3510 & 3774 & 4854 & 6018 & 6666 & 7470 & 10638 & 12534\\
13154 & 13206 & 20406 & 24270 & 32670 & 42186 & 45258 & 55878 & & \\
\hline 
\end{tabular}
\end{center}
\vspace{.2cm}

\centerline{$n \equiv 2$ (mod 12)}
\vspace{-.3cm}
\begin{center}
\begin{tabular}{|llllllllll|}
\hline 
13154 & & &&&&&&& \\
\hline
\end{tabular}
\end{center}
\caption{Missing integers  in the packing $(-1, 2, 2, 3)$ up to $10^6$ }
\end{table}

Our second example is the packing with root quadruple $(0,0,1,1)$.
As mentioned above, there are congruence conditions (mod $24$) in
this case. Table 4 presents numerical data on exceptional values
for the allowed congruence classes (mod $24$) up to $T = 10^7$.
There is a much larger set of exceptional values, and it appears
more equivocal whether the full list of exceptional values is
finite. However we think it is.

\begin{table}
\centerline{$n \equiv 0$ (mod 24)}
\begin{center}
\begin{tabular}{|lllllllll|}
\hline
 48    & 120   &  360  &  528 &  552&  720&  888 & 912 & 1080  \\
 1176  & 1272  &1392   &1560  &1704 & 1848  &1968  &2184  &2208 \\
 2736  & 2880  &3240   &3408  &3552 & 4080  &4392  &4464  &4584 \\
 4680  & 4896  &5040   &5088  &5760 & 6192  &6888  &7272  &8280 \\
 8880  & 9792  &10680  &10920  &10944  &11760  &11928  &13152  &14160 \\
 14328 & 16008 & 17160 & 17232 & 17520 & 18000 & 19320  &20712 & 23160 \\
 25896 & 26472 & 26760 & 27552 & 27600 & 27768 & 29424 & 29688  &30288 \\
 31440 & 34440 & 34488 & 35232 & 36408 & 36648 & 36816  &37968  &38928 \\
 39168 & 43056 & 43392 & 45240 & 46056 & 50448 & 52800  &58728  &59400 \\
 66120 & 74976 & 80280 & 82200 & 87192 & 93216  &96912  &96960  &107016\\
 108240& 117480& 121680 &133392 & 137280 & 138360 & 165360 & 201480&
 399000 \\
 424560 & 496080& & & & & & & \\ 
\hline 
\end{tabular}
\end{center}
\vspace{.2cm}

\centerline{$n \equiv 12$ (mod 24)}
\begin{center}
\vspace{-.3cm}
\begin{tabular}{|lllllllll|}
\hline
132&  252&  300&  468&  636&  780&  1140&  1476&  1572\\
1980&  2100&  2148&  2628&  2820&  2868&  3012&  3492&  3828\\
3900&  4212&  4692&  5028&  5148&  5340&  5796&  6516&  6684\\
6900&  7380&  7908&  8772&  10020&  10212&  10260&  10380&  10548\\
11268&  11868&  12876&  13572&  14100&  14244&  14724&  14916&
15300\\
15588&  19260&  19620&  20940&  21732&  22908&  23652&  24252&
24804\\
25140&  25812&  26100&  26124&  27660&  28860&  29532&  30540&
31092\\
31932&  36564&  37908&  38772&  39780&  41460&  41964&  44988&
46980\\
52260&  52788&  61596&  67308&  69324&  69420&  75900&  76908&
79740\\
88140&  101940&  120300&  135252&  185580&  188748&  220308&  228780&  
234660\\
354540&  422820&  472548&  926820&  1199820& & & & \\
\hline 
\end{tabular}
\end{center}
\vspace{.2cm}

\centerline{$n \equiv 1, 4$ or $9$ (mod 24)}
\vspace{-.3cm}
\begin{center}
\begin{tabular}{|l|llllll|ll|}
\hline 
241 & 340&  748&  2980&  5452&  11380&  45652& 16617&  21825\\
\hline
\end{tabular}
\end{center}
\vspace{.2cm}

\centerline{$n \equiv 16$ (mod 24)}
\vspace{-.3cm}
\begin{center}
\begin{tabular}{|llllllllll|}
\hline 
 208&  328&  712&  1168&  2488&  3400&  5200&  13600&  15088&  116896\\
\hline
\end{tabular}
\end{center}
\caption{Missing integers in the packing $(0, 0, 1, 1)$, up to $10^7$ }
\end{table}

The numerical examples above support the idea
that for any fixed integer
Apollonian packing and for sufficiently large integers a
finite list of 
congruence conditions will be  the only obstruction to
existence. We therefore propose the following strengthening
of the Density Conjecture.

\noindent {\bf Strong Density Conjecture}.
{\em In any primitive integral Apollonian packing, all sufficiently large
integers occur, provided they are not excluded by congruence conditions.
}

In further support of the Strong Density Conjecture, we note 
an analogy to
a number-theoretic conjecture of Zaremba \cite{zaremba}, who conjectured that
there exists an absolute constant $b$
(possibly $b=5$) such that each sufficiently large
positive integer can be represented by some 
continuant with digits bounded above
by $b$. In other words, given any integer $m>1$, 
there exists an integer $a<m$ ($a$ relatively prime to $m$)
such that the simple continued fraction $[0,c\sb 1,\cdots,c\sb r]=a/m$ 
has partial denominators $c\sb i \leq b$.
%Zaremba [Monatsh. Math. 78 (1974), 446--460;
Fix $b$, and  let $M$ be the set of all pairs $(a, m)$ with the above property.
There is a  linear recurrence for the pairs $(a, m)$ which is similar to
that of the Descartes quadruples, since if 
the terms in the
continued fraction of $a/m$ are bounded by $b$, 
then so are those for  the fractions $1/(i+a/m)$,
$i=1,2,\dots, b$. Zaremba's conjecture is saying that
all the integers $m$ will appear in some pair of $M$. 
This conjecture is currently still open. But as in the Apollonian
packing, consideration of the Hausdorff dimension 
of the set $E_b=\{ a/m: (a,m)\in M\}$ is suggestive. 
Namely, let $S_b(m)$ be the number of $a$'s such that $(a,m) \in M$. 
If $S_b(m) \sim m^\beta$, then
$\sum m^\beta m^{-x}$ converges iff  $x \geq \beta+1$.
Since the abscissa of convergence of the series 
$\sum S_b(m) m^{-x}$ is equal to twice the Hausdorff dimension $\gamma$
of $E_b$ (see T. Cusick \cite{Cu77}), then $\beta=2 \gamma -1
\approx .0624 >0$. Thus the ``expected'' number of appearances
of $m$ in the pairs of $M$ is $m^\beta \gg 1$. 

%--------------------------------------------------------------------
%
% Section 7
%
%--------------------------------------------------------------------
%

\section{The Growth of Descartes Quadruples in a Packing}
\setcounter{equation}{0}

The circles in an integral Apollonian circle packing, starting from
the root quadruple, are enumerated
by the elements of the Apollonian group. The graph of this group
is a rooted  infinite tree with four edges meeting each vertex,
with each vertex 
labelled by a nontrivial word in the generators of the Apollonian group.
(Such a word satisfies the condition that any two adjacent generators 
in the word are unequal.)
Starting from the root node, there are $4$ nodes at depth $1$, and at
each subsequent level there are three choices of generators at
each node, so there are
$4 \times 3^{n -1}$ words of length $n$ labelling depth $n$
circles. How are the curvatures of the circles at depth $n$ 
distributed? We consider the maximum value, the minimum value,
and the median value. In the process we also determine the 
joint spectral radius
of the generators of the Apollonian group.

We begin with the maximum value. We define for $n = 4m+i$ with 
$0 \le i \le 3$, the reduced word $T_n$ of length $n$ given by 
\beql{701a}
T_n := T_i (S_4 S_3 S_2 S_1 )^m,
\eeq
with $T_i = I,S_1, S_2 S_1, S_3 S_2 S_1$ for $0 \le i \le 3$, respectively.
\begin{theorem}\label{th71b}
Let $\bv = (a,b,c,d)$ be any root quadruple with 
$a \le b \le c \le d$ and $a < 0$, $a +b+c+d> 0$.
Then for any reduced word $W$ of length $n$ in the generators
$\{S_1, S_2, S_3, S_4 \}$ of the Apollonian group,
\beql{DI71}
|| W \bv ||_\infty \le ||T_n \bv ||_\infty ~.
\eeq
\end{theorem}

\pf Write
$W = S_{i_n} S_{i_{n-1}} \cdots S_{i_1}$ and set
$\bw^{(n)} = W \bv$ and $\bv^{(n)} =T_n \bv$.
Write the elements of $\bw^{(n)}$ and $\bv^{(n)}$ in increasing
order as
$$w_1^{(n)} \le w_2^{(n)} \le w_3^{(n)} \le w_4^{(n)} \quad\mbox{and}\quad
v_1^{(n)} \le v_2^{(n)} \le v_3^{(n)} \le v_4^{(n)} ~.
$$
The idea of the proof is that $T_n$ always inverts with respect to the 
circle of smallest curvature, and in fact produces
the largest curvature vector in a strong lexicographic sense.
More precisely, we prove by induction on $n \ge 1$ that
\beql{DI72}
w_i^{(n)} \le v_i^{(n)} \quad\mbox{for}\quad 1 \le i \le 4
\eeq
and
\beql{DI73}
w_4^{(n)} - w_1^{(n)} \le v_4^{(n)} - v_1^{(n)} ~.
\eeq
For the base case $n=1$, we have
$\bv^{(1)} = (a', b,c,d)$ where
$a' = 2(b+c+d) - a = ||S_1 \bv ||_\infty$.
If $b' = 2(a+c+d) -b = ||S_2 \bv ||_\infty$ and
$c' = || S_3 \bv ||_\infty$, $d' = ||S_4 \bv ||_\infty$ 
then $a \le b \le c \le d$
gives $d' \le c' \le b' \le a'$, and \eqn{DI72} holds for $n=1$,
since $d' \geq d$ because $\bv$ is a root quadruple. 

For the induction step, a reduced word has $i_n \neq i_{n-1}$.
The forbidden move $S_{i_{n-1}}$ is the one that replaces $w_4^{(n-1)}$ 
with $2(w_1^{(n-1)} + w_2^{(n-1)} + w_3^{(n-1)} ) - w_4^{(n-1)}.$
Now the induction hypothesis gives
\begin{eqnarray*}
w_1^{(n)} & \le & w_2^{(n-1)} \le v_2^{(n-1)} = v_1^{(n)} \\
w_2^{(n)} & \le & w_3^{(n-1)} \le v_3^{(n-1)} = v_2^{(n)} \\
w_3^{(n)} & \le & w_4^{(n-1)} \le v_4^{(n-1)} = v_3^{(n)}
\end{eqnarray*}
and
\begin{eqnarray*}
w_4^{(n)} & \le & 2 (w_2^{(n-1)} + w_3^{(n-1)} + w_4^{(n-1)} ) - w_1^{(n-1)} \\
& \le & 2(w_2^{(n-1)} + w_3^{(n-1)} ) + w_4^{(n-1)} + (w_4^{(n-1)} - w_1^{(n-1)} ) \\
& \le & 2(v_2^{(n-1)} + v_3^{(n-1)} ) + v_4^{(n-1)} + (v_4^{(n-1)} - v_1^{(n-1)} ) \\
& = & v_4^{(n)} ~.
\end{eqnarray*}
For the remaining inequality, suppose first that $w_1^{(n)} = w_2^{(n-1)}$.
Then
\begin{eqnarray*}
w_4^{(n)} - w_1^{(n)} & = &
[2(w_2^{(n-1)} + w_3^{(n-1)} + w_4^{(n-1)} ) - w_1^{(n-1)}] - w_2^{(n-1)} \\
& = &
w_2^{(n-1)} + 2w_3^{(n-1)} + w_4^{(n-1)} + (w_4^{(n-1)} -w_1^{(n-1)} ) \\
& \le &
v_2^{(n-1)} + 2v_3^{(n-1)} + v_4^{(n-1)} + (v_4^{(n-1)} -v_1^{(n-1)} ) \\
& = & v_4^{(n)} -v_1^{(n)} ~.
\end{eqnarray*}
If, however, $w_1^{(n)} = w_1^{(n-1)}$, then
\begin{eqnarray*}
w_4^{(n)} - w_1^{(n)} & \le &
[2(w_1^{(n-1)} + w_3^{(n-1)} + w_4^{(n-1)}) - w_2^{(n-1)} ] - w_1^{(n-1)} \\
& \le & 2(w_2^{(n-1)} + w_3^{(n-1)} + w_4^{(n-1)} ) - w_1^{(n-1)} - w_2^{(n-1)} \\
& \le &
v_4^{(n)} - v_1^{(n)} ~,
\end{eqnarray*}
using the previous inequality.
This completes the induction step.~~~$\Box$

The maximum growth rate of the elements at level $n$ of the Apollonian group is also describable in terms of the joint spectral radius of the generators
$\{S_1, S_2, S_3, S_4 \}$ of the Apollonian group.

\begin{defi}
\label{def71}
{\rm
Given a finite set of $n \times n$ matrices
$\Sigma = \{ M_1, \ldots, M_s \}$ the {\em joint spectral radius}
$\sigma ( \Sigma )$ is
$$\sigma ( \Sigma ) :=
\limsup_{k \to \infty} \left\{
\max_{1 \le i_1, \cdots, i_k \le s} \sigma (M_{i_1} \cdots M_{i_k} )^{1/k} \right\} ~,
$$
where 
$\sigma (M) := \max \{| \lambda | : \lambda ~\mbox{eigenvalue of $M$}\}$ 
is the spectral radius of $M$.
}
\end{defi}

The notion of joint spectral radius has appeared in many contexts, including 
wavelets and fractals; see Daubechies and Lagarias \cite{DL92} for a 
discussion and references. In general it is hard to compute, but here
we can obtain an explicit answer.

\begin{theorem}\label{th72b}
The joint spectral radius for the generators
$\Sigma = \{S_1, S_2, S_3, S_4 \}$ of the Apollonian group is
$\sigma (\Sigma ) = \theta^{1/4}$ where
\beql{DI74}
\theta = \frac{1}{2} \left( 1 + \sqrt{5} + \sqrt{2+2\sqrt{5}}\right) 
\approx 2.890 ~.
\eeq
It is attained by $M = S_4 S_3 S_2 S_1$.
\end{theorem}

\pf
Pick a fixed root quadruple with $a < 0$, say $\bv = (-1,2,2,3)$, and
 consider the associated packing $\sP_{\bv}.$
Lemma \ref{le52a} asserts that
\beql{DI75}
c_0 \| M \|_F \le || M \bv ||_\infty \le c_1 \|M \|_F , \quad\mbox{all}\quad
M \in \sA ~.
\eeq
We use the well-known fact that, for any real $n \times n$ matrix $M$,
\beql{DI76}
\sigma (M) = \lim_{k \to \infty} \| M^k \|_F^{1/k} ~.
\eeq
Now \eqn{DI75} gives for any reduced word
$M = S_{i_s} \cdots S_{i_2}S_{i_1} \in \sA$ with
$i_k \neq i_{k-1}$ that
$$\sigma (M)^{1/s} = \lim_{k \to \infty}
( || M^k \bv ||_\infty )^{\frac{1}{ks}} ~,
$$
Choosing $k=4n$, Theorem \ref{th71b} yields
$$\sigma (M)^{1/s} \le \lim_{n \to \infty} 
||T_{4ns} \bv ||_\infty^{\frac{1}{4ns}} ~.
$$
Since $T_{4ns} = (S_4 S_3 S_2 S_1 )^{ns}$, this gives
\begin{eqnarray*}
\sigma (M) & \le & \lim_{n \to \infty} 
||(S_4 S_3 S_2 S_1)^{ns} \bv ||_\infty^{\frac{1}{4ns}} \\
& \le & \sigma (S_4 S_3 S_2 S_1)^{1/4} ~.
\end{eqnarray*}
Choosing $M=S_4 S_3 S_2 S_1 \in \sA$ attains equality (with $s=4$), 
which determines the joint spectral radius.
A computation reveals that the characteristic polynomial of 
$M = S_4 S_3 S_2 S_1$ is $X^4 - 2X^3 - 2X^2 -2X +1 =0$ which factors as
$$(X^2 + (-1 + \sqrt{5} ) X +1) (X^2 - (1+ \sqrt{5} ) X+1) =0 ~.$$
Its spectral radius is given by \eqn{DI74}.~~~$\Box$

The minimal growth rate of any reduced word of length $2n$ is attained
 by the word
$W_{2n} = (S_4 S_3 )^n$.
If $\bv = (a,b,c,d)$ is a root quadruple with $a < 0$, then
the supremum norm
\beql{DI76a}
|| W_{2n} \bv ||_\infty = n(n+1) (a+b) - nc + (n-1) d
\eeq
grows quadratically with $n$.
We omit the easy proof of this fact.

To conclude this section, we consider the ``average value'' of 
$||W \bv ||_\infty$ over all reduced words
$W$ in $\sA$ of length $n$, which we define to
 be the {\em median} of this distribution. (The elements of the
distribution are exponentially large, so the median is a more appropriate
quantity to consider than the mean value.)
Let $T_n$ denote the median.
We expect that its growth rate  should be related
 to the Hausdorff dimension $\alpha$ of the 
limit set of the Apollonian packing.
The results of \S5 lead to the heuristic that one should expect
\beql{DI77}
\lim_{n \to \infty} \frac{1}{n} \log T_n = \frac{\log 3}{\alpha} ~.
\eeq
We leave the proof (or disproof) of this as an open problem.
 
%
% Section 8
%
 
\section{Open questions}
\setcounter{equation}{0}
There remain many open  questions concerning
integral Apollonian circle packings. 
We list a few of these here.

(1) In any primitive integral
Apollonian packing $\mathcal P$, just four distinct residue
classes modulo $12$ can occur as curvature values in $\mathcal P$.
For example, for $\mathcal P = (0,0,1,1)$, these values are
$0,1,4,9~(\bmod~12)$ while for $\mathcal P = (-1,2,2,3)$, 
they are $2,3,6,11~(\bmod~12)$. As we noted in Section $6$, it seems
likely that in the packing $(-1,2,2,3)$, all sufficiently large
integers congruent to $2,3,6$ and $11~(\bmod~12)$ actually do occur.
However, in the packing $(0,0,1,1)$, instead of the $8$ residue
classes $0,1,4,9,12,13,16,22~(\bmod~24)$ which we might expect
to occur, the classes $13$ and $22$ are completely missing.
Again, computation suggests that only finitely many values in the
other $6$ are missing in $(0,0,1,1)$. Is it true that in any
integral Apollonian packing, the only congruence restrictions
 on the curvature values are for the modulus 24?
However, in no case can we even show that the set of values
which do occur has positive upper density.

(2) Is there a direct way for 
determining the root quadruple to which
a given Descartes quadruple belongs? The only way we currently know 
uses the reduction algorithm described in Section $3$.
 
(3) Concerning root quadruples, what are the asymptotics of the
total number of root quadruples having Euclidean height below $T$?

(4) We have not
proved any reasonable lower bound on the number of integers
below $T$ that occur as curvatures in a fixed integral Apollonian
packing. For how large a $\beta$ can one  prove asymptotically
that
at least $T^{\beta + o(1)}$ integers occur in every such packing?

(5) All of the preceding questions can also be raised for
integral Apollonian packing of spheres in $3$ dimensions,
as discussed in \cite{GLMWY23}.
 For example,
what are the modular restrictions (if any) for the Descartes quintuples
$(a,b,c,d,e)$ occurring in the packing with root quintuple $(0,0,1,1,1)$?

(6) There exist strongly integral
Apollonian packings, in which the circles all have integer curvatures
and also the   curvature$\times$centers of the circles are Gaussian integers,
where the circle centers are coordinatized as complex numbers. 
(See \cite{GLMWY22}.)
Most  questions investigated for integral packings
can also be asked for strongly integral packings.
If we write $(x, xX)$ for
a circle with curvature $x$ and (complex) center $X$, then
the pairs $(x, xX)$ must also satisfy various modular constraints.
For example, modulo $12$, the standard integral packing (i.e., 
starting with the circles $(-1,0), (2,1), (2,-1)$ has just
$20$ types of circles, namely,
\begin{eqnarray*}
(x,xX) &=&
(2,1), (2,3), (2,5), (2,7), (2,9), (2,11)\\
&&(3,2i), (3,4i), (3,8i), (3,10i)\\
&&(6,3+4i), (6,3+8i), (6,9+4i), (6,9+8i)\\
&&(11,0), (11,4), (11,8), (11,6i), (11,4+6i), (11,8+6i)
\end{eqnarray*}
and there are just $120$ different four-circle configurations.
What are the asymptotics of these types and configurations?
What is the characterization of the integral (complex) vectors
$(x, xX)$ that can appear in a given packing?

We hope to return to some of these issues in a future paper.

\paragraph {Acknowledgments.} The authors wish to acknowledge the 
insightful comments of 
Arthur Baragar, William Duke, Andrew Odlyzko, Eric Rains, and 
Neil Sloane at various stages of this work. We also thank 
the referee for many useful comments and historical references.

{\tt
\begin{tabular}{lllll}
email: & graham@ucsd.edu \\
& jcl@research.att.com \\
  & colinm@research.avayalabs.com   \\
 & allan@research.att.com \\
& Catherine.Yan@math.tamu.edu 
\end{tabular}
 }
\end{document}